\documentclass[a4paper,10pt]{article}

\usepackage{amssymb,latexsym}
\usepackage{euscript}
\usepackage{eufrak}
\usepackage{amsmath}
\usepackage{tikz}
\usepackage{bussproofs}
\usepackage{mathrsfs} 

\usetikzlibrary{matrix}

\newcommand{\mm}{\mathcal} 
\newcommand{\rr}{\mathscr} 
\newcommand{\bb}{\mathbb}  
\newcommand{\ff}{\mathfrak} 

\newcommand{\pt}{\forall} 
\newcommand{\ex}{\exists} 
\newcommand{\no}{\neg} 
\newcommand{\imp}{\rightarrow} 
\newcommand{\sii}{\leftrightarrow} 

\newcommand{\und}{\frac 1 2} 

\newcommand{\sneg}{{\sim}} 
\newcommand{\cons}{{\circ}} 

\newcommand{\lfdash}{\Vdash_{\rm LFI1}}

\newcounter{defcounter}
\newcommand{\uno}[1]{#1\wedge\cons #1}
\newcommand{\cero}[1]{\no #1\wedge\cons #1}
\newcommand{\tres}[1]{#1\wedge\no #1}

\newcommand{\mas}[1]{\|#1\|^{\mathfrak A}_\oplus}
\newcommand{\por}[1]{\|#1\|^{\mathfrak A}_\odot}
\newcommand{\menos}[1]{\|#1\|^{\mathfrak A}_\ominus}

\newcommand{\mass}[2]{\|#1\|^{\mathfrak #2}_\oplus}
\newcommand{\porr}[2]{\|#1\|^{\mathfrak #2}_\odot}
\newcommand{\menoss}[2]{\|#1\|^{\mathfrak #2}_\ominus}

\newcommand{\tas}[1]{#1^{\mathfrak A}_\oplus}
\newcommand{\tpor}[1]{#1^{\mathfrak A}_\odot}
\newcommand{\tenos}[1]{#1^{\mathfrak A}_\ominus}

\newcommand{\triple}[1]{#1^{\mathfrak A}=\langle #1^{\mathfrak A}_\oplus, #1^{\mathfrak A}_\ominus, #1^{\mathfrak A}_\odot \rangle}

\newcommand{\triplee}[2]{#1^{\mathfrak #2}=\langle #1^{\mathfrak #2}_\oplus, #1^{\mathfrak #2}_\ominus, #1^{\mathfrak #2}_\odot \rangle}

\newcommand{\ter}[1]{t_{ #1}^\ff A[s]}

\newcommand{\sent}[1]{{\rm Sent}(#1)}
\newcommand{\eledos}{{\bf LFI2}}
\newcommand{\tripc}[1]{\|#1\|^{\mathfrak A}}

\newcommand{\dodash}{\Vdash_\quci}

\newcommand{\vct}[1]{#1_1,\ldots ,#1_n }

\newcommand{\acio}[1]{{\bf (Ax#1)}}

\newcommand{\brac}[1]{\langle #1\rangle}

\newcommand{\mapa}[3]{#1:#2\rightarrow #3}
\newcommand{\embe}[4]{#1:#2\prec_{#4}#3}

\newcommand{\sisi }[1]{\{#1,\cons #1\}}
\newcommand{\nosi}[1]{\{#1,\no #1\}}
\newcommand{\nono}[1]{\{\no #1,\cons #1\}}

\newenvironment{myequationbc1}{%

\begin{equation}}
{\end{equation}}

\newenvironment{myequationci}{%

\begin{equation}}
{\end{equation}}

\newenvironment{myequationcf}{%

\begin{equation}}
{\end{equation}}

\newenvironment{myequationce}{%

\begin{equation}}
{\end{equation}}

\newenvironment{myequationco1}{%

\begin{equation}}
{\end{equation}}
\newenvironment{myequationco2}{%

\begin{equation}}
{\end{equation}}
\newenvironment{myequationco3}{%

\begin{equation}}
{\end{equation}}

\newenvironment{myequationcr1}{%

\begin{equation}}
{\end{equation}}
\newenvironment{myequationcr2}{%

\begin{equation}}
{\end{equation}}
\newenvironment{myequationcr3}{%

\begin{equation}}
{\end{equation}}

\newenvironment{myequation}{%
\addtocounter{equation}{-1}
\refstepcounter{defcounter}

\begin{equation}}
{\end{equation}}

\newcommand{\ese}[1]{S(\ff #1)}

\newtheorem{definition}{\vspace{1mm}Definition}[section]
\newtheorem{remark}[definition]{\vspace{1mm}Remark}
\newtheorem{lemma}[definition]{\vspace{1mm}Lemma}
\newtheorem{proposition}[definition]{\vspace{1mm}Proposition}
\newtheorem{theorem}[definition]{\vspace{1mm}Theorem}
\newtheorem{corollary}[definition]{\vspace{1mm}Corollary}

\newenvironment{dem}{\noindent\bf Proof. \rm}{\hfill $\mbox{\boldmath{$ \square$}}$}

\newcommand{\ciore}{{\bf Ciore}}
\newcommand{\quci}{{\bf QCiore}}
\newcommand{\lfis}{{\bf LFI}s}
\newcommand{\lfi}{{\bf LFI}}
\newcommand{\mbc}{{\bf mbC}}
\newcommand{\qmbc}{{\bf QmbC}}
\newcommand{\dacdot}{{\bf J3}}

\makeatother

\makeatletter

\newcommand*{\reflectmathsymbol}[2]{%
  \reflectbox{$\m@th#1#2$}%
}

\makeatletter
\providecommand*{\Dashv}{%
  \mathrel{%
    \mathpalette\@Dashv\models
  }%
}
\newcommand*{\@Dashv}[2]{%
  \reflectbox{$\m@th#1#2$}%
}

\makeatother

\title{Some model-theoretic results on the 3-valued\\ paraconsistent first-order logic {\bf QCiore}}

\author{
    Marcelo E. Coniglio $^\textup{\scriptsize a}$
    \and Tadeo G. Gomez $^\textup{\scriptsize b}$
    \and Mart\'{\i}n Figallo  $^\textup{\scriptsize b}$
}

\date{
    $^\textup{\scriptsize a}$\textit{\small Institute of Philosophy and the Humanities (IFCH), and\\
\small Centre for Logic, Epistemology and the History of Science (CLE)\\
\small University of Campinas (UNICAMP), Brazil\\
\small {\em Email: \texttt{coniglio@cle.unicamp.br}}}
    \\
    $^\textup{\scriptsize b}$\textit{\small Departament of Mathematics and Instituto de Matem\'atica  de Bah\'ia Blanca (INMABB)\\  National University of the South (UNS), Bah\'{\i}a Blanca, Argentina\\
\small {\em Email: \texttt{$\{$tadeogerman, figallomartin$\}$@gmail.com}}}
}

\begin{document}

\maketitle

\begin{abstract}
In this paper the 3-valued paraconsistent first-order logic \quci\ is studied from the point of view of Model Theory. The  semantics for \quci\ is given by partial structures. They are first-order structures such that each $n$-ary predicate $R$ is interpreted as a triple of pairwise disjoint sets of $n$-uples representing, respectively, the set of tuples which  actually belong to $R$,  the set of tuples which actually do not belong to $R$, and the set of tuples whose status is dubious or contradictory. Partial structures were proposed in 1986 by I. Mikenberg, N. da Costa and R. Chuaqui for the theory of quasi-truth (or pragmatic truth). In 2014, partial structures were studied by M. Coniglio and L. Silvestrini  for a 3-valued paraconsistent first-order logic called  {\bf LPT1}, whose 3-valued propositional fragment is equivalent to da Costa-D'Otaviano's logic \dacdot.  This approach is adapted in this paper  to   \quci, and some important results of classical Model Theory such as Robinson's joint consistency theorem, amalgamation and interpolation are obtained.  Although we focus on \quci, this framework can be adapted to other 3-valued first-order logics.

\vspace*{4mm}


\vspace*{4mm}

\end{abstract}

\section{\large \bf Introduction}\label{s1}

A logic is paraconsistent if it has a negation which allows the existence of contradictory but nontrivial theories. The first sistematic study of paraconsistent logics was given by N. da Costa in 1963, when presented in~\cite{daC} his well-known hierarchy $C_n$ (for $n \geq 1$) of $C$-systems. His approach to paraconsistency is based on the idea that propositions are {\em prima facie} {\em dubious}, in the sense that a sentence $P$ and its negation $\neg P$ can be assumed to hold simultaneously, without trivialization.  The sentence $P^\circ:=\neg (P \wedge \neg P)$ indicates the {\em well-behavior} (in da Costa's terminology) of a sentence $P$ in the stronger system $C_1$.\footnote{For each $n \geq 1$ there is a sentence $P^{\circ_n}$ denoting the well-behavior of $P$.} That sentence is used  to express that a sentence $P$ is {\em reliable} (meaning that either $P$ or $\neg P$ cannot be negated withou trivializing). This being so, from $P$, $\neg P$ and $P^\circ$ (as well as from  from $\neg P$, $\neg\neg P$ and $P^\circ$) every sentence $Q$ can be derived. In other words: a theory contaning a contradiction $\{P,\neg P\}$ is not always deductively trivial (since $P$ represents, {\em prima facie}, dubious information, hence it could be contradictory) but any theory containing  $\{P,\neg P, P^\circ\}$ is always deductively trivial (since it contains a contradiction involving absolutely reliable information, which is an absurd). 

The approach to paraconsistency of da Costa, nowadays known as {\em Brazilian school of paraconsistency}, was naturally generalized by W. Carnielli and J. Marcos with the notion of  {\em Logics of Formal Inconsistency} ({\bf LFI}s for short). These logics are paraconsistent and, additionally, have a unary connective $\circ$ (primitive or definable) that allows us
to express the notion of consistency of sentences inside the object language. Thus, $\circ P$ means that {\em $P$  is consistent}. The  {\em  Explosion Law} (`from $P$ and $\neg P$ every $Q$ follows') does not hold
in \lfis, as much as in any other paraconsistent logic. But  some contradictions will cause deductive explosion: contradictions involving a  {\em  consistent} sentence lead to triviality -- intuitively, and as in da Costa's $C$-systems, one can understand this situation as a contradiction involving well-established facts,
or involving propositions that have conclusive favorable evidence (what da Costa called {\em well-behaved} sentences). In this
sense, \lfis\ -- including da Costa's systems --  are logics that allow us to separate the sentences for which the Explosion Law hold, from those for which  does not hold. The only difference between da Costa's $C$-systems and \lfis\ in general is that, in the latter, the consistency operator can be possibly given by an unary connective $\circ$ while, in the former, it is always defined in terms of the other connectives.

In \cite{it127, dot:85a, dot:85b, dot:87}, I. D'Ottaviano developed a 3-valued model theory in order to study a first order version of the logic \dacdot, introduced by herself and N. da Costa in~\cite{it131}. The logic \dacdot\ is a 3-valued paraconsistent propositional logic; moreover, it is an {\bf LFI}, as it can easily proven (see~\cite{Tax}). In~\cite{lf79}, W.~Carnielli, J.~Marcos  and S.~de~Amo propose the use of \dacdot, seen as an \lfi\ under the name of {\bf LFI1}, in order to deal with inconsistent databases. Besides this, they define a first-order version of {\bf LFI1} called ${\bf LFI1^*}$, which happens to be equivalent to the first-order version of \dacdot\ studied by D'Ottaviano, and it is also equivalent to M. Coniglio and L. Silvestrini's logic {\bf LPT1} (see~\cite{consil}). In the same paper it was studied  another 3-valued \lfi\ called {\bf LFI2} (renamed as \ciore\ in~\cite{Tax}), as well as a first-order version of it called ${\bf LFI2^*}$. As it will be shown in Remark~\ref{remquanti}, the first-order version ${\bf LFI2^*}$ of \ciore\ proposed in~\cite{lf79} differs from the one (\quci) proposed in the present paper.
 
The logic \ciore\ is an axiomatic extension of \mbc, the basic \lfi\ considered in~\cite{WCMCJM} and~\cite{CarCon}, with strong properties concerning propagation of the consistency operator $\circ$. In addition to its syntactical presentation by means of a Hilbert calculus, it was proved in~\cite{lf79} that \ciore/{\bf LFI2} is sound and complete with respect to a 3-valued matrix.

Positive classical logic ({\bf CPL}$^+$) has been assumed as a natural starting point from which many {\bf LFI}s are defined, although 
some {\bf LFI}s are studied starting from other logics than {\bf CPL}$^+$. In this sense, {\bf mbC} is the basic logic of formal inconsistency that, starting with positive classical logic {\bf CPL}$^+$ and adding a negation and a consistency operator, it is endowed with minimal properties in order to satisfy the definition of {\bf LFI}s (see~\cite{WCMCJM} and~\cite[Chapter~2]{CarCon}).

Besides being a subsystem of {\bf CPL}, {\bf mbC} is also an extension of
{\bf CPL} obtained by adding to the latter a consistency operator $\cons$  and a paraconsistent
negation $\neg$. In this sense, {\bf mbC} can be viewed, both, as a subsystem
and as a conservative extension of {\bf CPL}. A similar phenomenon holds for several
other {\bf LFI}s. A remarkable feature of {\bf LFI}s in general, and of {\bf mbC} in particular, is that classical logic ({\bf CPL}) can be codified, or recovered, inside such logics. A first-order extension of \mbc\ called \qmbc\ was proposed in~\cite{car:etal:14}, and it was also studied in~\cite[Chapter~7]{CarCon}. In~\cite{fer:18}, T. Ferguson has obtained an important  model-theoretic result on \qmbc,  the Keisler-Shelah Theorem, by using a novel technique of model-theoretic ``atomization''.

The question of characterizability (or not) by finite matrices, as well as the
algebraizability of (extensions of) {\bf mbC} has been tackled in~\cite{WCMCJM} and~ \cite{CarCon}. Some negative results, in the
style of the famous Dugundji's theorem for modal logics, were shown for several
extensions of {\bf mbC}. These results have shown that a wide variety of {\bf LFI}s extending {\bf mbC} cannot be
semantically characterized by finite matrices. Despite these general result, there are many
 {\bf LFI}s which, besides being characterizable by a single 3-valued matrix, are algebraizable in the well-known sense of Blok and Pigozzi. 
Among them, we can mention
Halld\'en's logic of nonsense as well as Segerberg's variation,  \dacdot, Sette's logic
{\bf P1} and the system \ciore, to be studied in this paper.

The  logic \dacdot\  has been re-introduced and studied independently by several authors at several times, and with different motivations. For instance, and as mentioned above, Carnielli, Marcos and de Amo introduced in~\cite{lf79} a propositional logic, which they called  {\bf LFI1}, as a suitable tool for modeling processes in databases with contradictory information. They proved the equivalence (up to language) between  \dacdot\ and {\bf LFI1}.
In \cite{cluns}, D. Batens and K. De Clercq studied the logic {\bf CluNs} which was introduced syntactically by means of a first order Hilbert calculus. The propositional fragment of {\bf CluNs} turned out to be equivalent to \dacdot. 
Finally, in \cite{consil}, Coniglio and Silvestrini presented the logic {\bf MPT} pursuing to establish a 3-valued matrix associated to pragmatic structures which constitute a generalization of the notion of {\em quasi-truth} introduced by Mikenberg et al. (see \cite{qt224}). Also, they prove that {\bf MPT} and \dacdot\ (and therefore also {\bf LFI1}) are functionally equivalent.
Originally, the  signature of {\bf LFI1} had the primitive connective  $\bullet$ (called {\em inconsistency operator}) to mark out inconsistencies in databases, by means of the equivalence 

$$\bullet \varphi\leftrightarrow (\varphi\wedge\neg \varphi)$$

\noindent In \cite[Chapter~4]{CarCon}, Carnielli and Coniglio introduced a new presentation {\bf LFI1'} of {\bf LFI1} by considering the connective $\cons$ (the consistency operator)  instead of $\bullet$ (the inconsistency operator), together with a suitable  Hilbert calculus called ${\bf LFI1}_\circ$. The logic {\bf LFI1'}   will be considered  in Remark~\ref{ejP1LFI1} below.

The semantic notion of quasi-truth, defined by means of partial structures, was  introduced in 1986 by I. Mikenberg, N. da Costa and R. Chuaqui (cf. \cite{qt224}). Partial structures are first-order structures where each $n$-ary predicate $R$ is interpreted as a triple  $\langle R_{+}, R_{-}, R_{u}\rangle$ of pairwise disjoint sets of $n$-uples over the domain of the structure such that $R_{+}$ is the set of tuples which  actually belong to $R$, $R_{-}$ is  the set of tuples which actually do not belong to $R$, and $R_{u}$ is the set of tuples whose status is dubious or contradictory.

In this paper a model-theoretic framework based on triples is  proposed for \quci, the first-order version of  \ciore. This approach is based on the pragmatic structures semantics proposed for {\bf LPT1} in~\cite{consil}.  Besides adequacy of  \quci\ w.r.t. pragmatic structures, some important results of classical Model Theory are obtained for this logic, such as Robinson's joint consistency theorem, amalgamation and interpolation.  Although we focus on \quci, this framework can be adapted to other 3-valued first-order {\bf LFI}s.

\

\section{Preliminaries}\label{s1.1}

\

In this section, a Hilbert-style presentation for {\bf Ciore} will be given. We consider a propositional language $\rr L$ over the propositional signature $\{\wedge,\vee,\imp,\no,\cons\}$ where the formulas are constructed as usual from a denumerable set $\{p_1,p_2, \ldots\}$ of propositional variables. 

The propositional logic \ciore \ is defined over the language $\rr L$ by the following Hilbert system:


\newpage 

\noindent{\bf Axiom schemata:} 

\vspace{-0.3cm}
\begin{myequation}
  \alpha\imp(\beta\imp\alpha)
   \end{myequation}
\vspace{-0.5cm}
\begin{myequation}
  \big(\alpha\imp(\beta\imp\gamma)\big)\imp\big((\alpha\imp\beta)\imp(\alpha\imp\gamma)\big)
   \end{myequation}
\vspace{-0.5cm}
\begin{myequation}
  \alpha\imp\big(\beta\imp(\alpha\wedge\beta)\big)
   \end{myequation}
\vspace{-0.5cm}
\begin{myequation}
  (\alpha\wedge\beta)\imp\alpha
   \end{myequation}
\vspace{-0.5cm}
\begin{myequation}
  (\alpha\wedge\beta)\imp\beta
   \end{myequation}
\vspace{-0.5cm}
\begin{myequation}
  \alpha\imp(\alpha\vee\beta)
   \end{myequation}
\vspace{-0.5cm}
\begin{myequation}
  \beta\imp(\alpha\vee\beta)
   \end{myequation}
\vspace{-0.5cm}
\begin{myequation}
  \big(
\alpha\imp\gamma
\big)
\imp
\big(
(\beta\imp\gamma)\imp((\alpha\vee\beta)\imp\gamma)\big
)
   \end{myequation}
\vspace{-0.5cm}
\begin{myequation}
  (\alpha\imp\beta)\vee\alpha
   \end{myequation}
\vspace{-0.5cm}
\begin{myequation}
  \alpha\vee\no\alpha
   \end{myequation}
\vspace{-0.5cm}
\begin{myequationbc1}
\cons\!\alpha\imp\big(\alpha\imp(\no\alpha\imp\beta)\big)
\end{myequationbc1}
\vspace{-0.5cm}
\begin{myequationci}
\no\cons\!\alpha\imp(\alpha\wedge\no\alpha)
\end{myequationci}
\vspace{-0.5cm}
\begin{myequationcf}
\no\no\alpha\imp\alpha
\end{myequationcf}
\vspace{-0.5cm}
\begin{myequationce}
\alpha\imp\no\no\alpha
\end{myequationce}
\vspace{-0.5cm}
\begin{myequationco1}
(\cons\alpha\vee\cons\beta)\imp\cons(\alpha\wedge\beta)
\end{myequationco1}
\vspace{-0.5cm}
\begin{myequationco2}
(\cons\alpha\vee\cons\beta)\imp\cons(\alpha\vee\beta)
\end{myequationco2}
\vspace{-0.5cm}
\begin{myequationco3}
(\cons\alpha\vee\cons\beta)\imp\cons(\alpha\imp\beta)
\end{myequationco3}
\vspace{-0.5cm}
\begin{myequationcr1}
\cons(\alpha\wedge\beta)\imp(\cons\alpha\vee\cons\beta)
\end{myequationcr1}
\vspace{-0.5cm}
\begin{myequationcr2}
\cons(\alpha\vee\beta)\imp(\cons\alpha\vee\cons\beta)
\end{myequationcr2}
\vspace{-0.5cm}
\begin{myequationcr3}
\cons(\alpha\imp\beta)\imp(\cons\alpha\vee\cons\beta)
\end{myequationcr3}

\noindent{\bf Inference rule:}
\begin{center}
{\bf (MP)}\hfill $\displaystyle\frac{\alpha,\alpha\imp\beta}{\beta}$\hfill\hspace{1mm}
\end{center}

\

\noindent Observe that axioms {\bf (Ax1)}-{\bf (Ax9)} plus {\bf (MP)} constitute a Hilbert calculus for positive classical logic ({\bf CPL}$^+$).
The logic \ciore\ was proposed as a suitable extension of
{\bf Cio} which is algebraizable in the sense of Blok and Pigozzi (see \cite{CarCon}).
Morever, as it was observed in \cite{Tax}, the logic \ciore  \, is so strong that it can be
characterized by a 3-valued matrix logic.

\begin{theorem}\label{SoundCompFI2}(See \cite[Theorem~4.4.29]{CarCon})
The system \ciore\ is sound and complete with respect to the following three valued matrix $M_e =\langle M, D \rangle$  over the signature $\Sigma$ with domain $M=\{1,\frac{1}{2},0\}$ and set of designated values $D=\{1,\frac{1}{2}\}$ such that
the truth-tables associated to each connective are the following:

\begin{center}
$
\begin{array}{|c||c|c|c|}
\hline
\wedge &1 &\frac{1}{2} &0 \\
\hline
\hline
1& 1 & 1 & 0 \\
\hline
\frac{1}{2} & 1 & \frac{1}{2} & 0 \\
\hline
0&0&0&0\\
\hline
\end{array}
$ \quad
$
\begin{array}{|c||c|c|c|}
\hline
\vee &1 &\frac{1}{2}&0 \\
\hline
\hline
1&1&1&1\\
\hline
\frac{1}{2}&1&\frac{1}{2}& 1 \\
\hline
0&1& 1 &0\\
\hline
\end{array}
$ \quad
$
\begin{array}{|c||c|c|c|}
\hline
\imp &1 &\frac{1}{2} &0 \\
\hline
\hline
1&1& 1 &0\\
\hline
\frac{1}{2}&1&\frac{1}{2}&0\\
\hline
0&1&1&1\\
\hline
\end{array}
$
\quad
$
\begin{array}{|c||c||c|}
\hline
 &\no &\cons \\
\hline
\hline
1&0&1 \\
\hline
\frac{1}{2}&\frac{1}{2}&0\\
\hline
0&1&1 \\
\hline
\end{array}
$
\end{center}

\end{theorem}

\

\

\noindent As mentioned in the Introduction, the logic associated with this matrix was proposed in~\cite{lf79} under the name {\bf LFI2}, being renamed as \ciore\ in~\cite{Tax}. In the following result, we exhibit a list of theorems of {\bf Ciore} that we shall use along this paper. We will write \ $\sneg\alpha$ \ as an abbreviation of \ $\no\alpha\wedge\cons\alpha$ to denote the strong (or classical) negation definable in \ciore.

\begin{theorem}\label{esquema} The following schemas are theorems of {\bf Ciore}. 

\hspace{-0.5cm}
\begin{tabular}{llll}
(i) & $\alpha\imp\alpha$, & (ii) & $(\alpha\imp\beta)\imp(\sneg\beta\imp\sneg\alpha)$, \\
(iii) & $(\alpha\wedge\neg\alpha)\sii\sneg{\cons}\alpha$, & (iv) & $(\alpha\wedge\neg\alpha)\sii\neg{\cons}\alpha$, \\
(v) & ${\cons}{\cons}\alpha$, & (vi) &  $\cons\alpha\sii\cons\neg\alpha$, \\
(vii) & $(\sneg\alpha \vee \sneg\beta) \sii
\sneg(\alpha\wedge\beta)$,  &
(viii) & $\big((\alpha\wedge\neg\alpha)\wedge(\beta\wedge\neg \beta)\big) 
\sii \big((\alpha\wedge\beta)\wedge\neg(\alpha\wedge\beta)\big)$,\\
\end{tabular}\\
\hspace*{5mm}
(ix) \ $\big((\alpha\wedge(\beta\wedge\cons\beta)\big) \vee ((\alpha\wedge\cons\alpha)\wedge\beta)\big)  \sii \big((\alpha\wedge\beta) \wedge \cons(\alpha\wedge\beta)\big)$,\\
(x) \ $\big((\alpha\wedge \cons \alpha) \vee(\beta\wedge\cons\beta) \vee \big((\alpha\wedge \neg \alpha) \wedge\sneg\beta\big) \vee \big(\sneg\alpha \wedge(\beta\wedge\neg\beta)\big)\big)  \sii \big((\alpha\vee\beta) \wedge \cons(\alpha\vee\beta)\big)$,\\
(xi) \ $\big(\sneg\alpha \wedge \sneg\beta\big)  \sii \sneg(\alpha\vee\beta)$,\\
(xii) \ $\big((\alpha\wedge \neg\alpha) \wedge(\beta\wedge\neg\beta)\big)  \sii \big((\alpha\vee\beta) \wedge \neg(\alpha\vee\beta)\big)$,\\
(xiii) \ $\big(\sneg\alpha \vee(\beta\wedge\cons\beta)\big)  \sii \big((\alpha\to\beta) \wedge \cons(\alpha\to\beta)\big)$,\\
(xiv) \ $\big(\alpha\wedge \sneg\beta\big)  \sii \sneg(\alpha\to\beta)$,\\
(xv) \ $\big((\alpha\wedge \neg\alpha) \wedge(\beta\wedge\neg\beta)\big)  \sii \big((\alpha\to\beta) \wedge \neg(\alpha\to\beta)\big)$.
\end{theorem}
\begin{dem} It is routine to check that all these formulas are tautologies in {\bf LFI2} by observing that, for every $x \in M$: $x = 1$ iff $x \wedge \cons x \in D$; $x = 0$ iff $\neg x \wedge \cons x \in D$; and $x = \frac{1}{2}$ iff $x \wedge \neg x \in D$. The result follows by Theorem~\ref{SoundCompFI2}.
\end{dem}

\

\section{Partial relations} \label{s1.5}

In order to define the sematical framework for \quci, the first-order version of \ciore\ proposed in this paper, it is necessary to consider partial relations (or triples).

Let ${\bf 3}=\{0,\frac{1}{2},1\}$ and consider the algebraic structure $\langle {\bf 3}; \vee, \wedge, \imp, \neg, \cons \rangle$ where the operations $\vee, \wedge, \imp,$ $\neg$ and $\cons$ are defined as in Theorem \ref{SoundCompFI2}. Let $X$ be a non-empty set. A {\em triple over $X$}  is a map $\mapa rX{{\bf 3}}$. If $r$ is a triple over $X$ we write
$r= \langle r_{ \oplus},r_{ \ominus},r_{ \odot}\rangle$
where  \  \
$r_{ \oplus}=r^{-1}(1)$, \ 
$r_{ \odot}=r^{-1}(\frac{1}{2})$, \
$r_{ \ominus}=r^{-1}(0)$. 

As usual, we denote by ${\bf 3}^X$ the set of all triples over $X$.

\begin{remark} \label{op3X}
Clearly, ${\bf 3}^X$ inherits the algebraic structure of ${\bf 3}$ where the  operations are defined pointwise. Namely, $(r \div u)(x) := r(x) \div u(x)$ and $(*r)(x):=*(r(x))$ for $\div \in \{\wedge,\vee,\to\}$, $* \in \{\neg,\cons\}$, $x \in X$ and $r,u \in  {\bf 3}^X$.
\end{remark}

\begin{proposition}\label{lfi1al}
Let $r = (r_\oplus,r_\ominus,r_\odot)$ and  $u = (u_\oplus,u_\ominus,u_\odot)$ two triples over $X$. Then:
\begin{itemize}
\item[(i)] $r\wedge u = ( ( r_\oplus \cap u_\oplus)\cup( r_\oplus\cap  u_\odot)\cup( r_\odot\cap  u_\oplus),  r_\ominus\cup  u_\ominus,  r_\odot \cap  u_\odot)$,
\item[(ii)] $ r \vee u  = (  r_\oplus\cup u_\oplus, r_\ominus\cap u_\ominus, (r_\odot\cap u_\ominus) \cup (r_\ominus\cap u_\odot) \cup (r_\odot\cap u_\odot))$,
\item[(iii)] $ r\imp u = ( r_\ominus\cup u_\oplus,  (r_\oplus\cup r_\odot)\cap u_\ominus, (r_\oplus\cup r_\odot)\cap u_\odot) $,
\item[(iv)] $\no r = (r_\ominus,r_\oplus,r_\odot)$,
\item[(v)] $\cons r = (r_\oplus \cup r_\ominus, r_\odot,\emptyset)$.
\end{itemize}
\end{proposition}
\begin{dem}  Here, only item~(ii) will be proved. First, $x\in (r\imp u)_\oplus$  iff  $r(x)\imp u(x)=1$  iff  $r(x)=0$ or $u(x)=1$   iff  $x\in r_\ominus\cup r_\oplus$. Secondly, $x\in (r\imp u)_\ominus$  iff  $r(x)\imp u(x)=0$  iff  $r(x)\in\{1,\und\}$ and $u(x)=0$ iff $x\in (r_\oplus\cup r_\odot)\cap u_\ominus$. Finally, $x\in (r\imp u)_\odot$  iff $r(x)\imp u(x)=\und$ iff $r(x)\in\{1,\und\}$ and $u(x)=\und$ iff $x\in (r_\oplus\cup r_\odot)\cap u_\oplus$.
\end{dem}

\

\noindent
The notion of {\em Partial Relation} introduced by Mikenberg et al. (see \cite{qt224}) is analogous to the notion of triple.
We interpret the nature of a triple $r = (r_\oplus,r_\ominus,r_\odot)$ over $X\neq\emptyset$ in the context of \quci\ as being a kind of 3-valued fuzzy set such that:
\ 
\begin{itemize}
\item[-] $r_{ \oplus}$ is the set of individuals which belongs to $r$,
\item[-] $r_{ \ominus}$ is the set of individuals which do not belongs to $r$,
\item[-] $r_{ \odot}$ is the set of individuals whose status w.r.t. $r$ is dubious or contradictory.
\end{itemize}

\noindent
It is not difficult to see that, if $r \in {\bf 3}^X$, then: (1)~$r_\oplus\cup r_\ominus\cup r_\odot=X$; and (2)~$r_\oplus\cap r_\ominus=r_\oplus\cap r_\odot=r_\ominus\cap r_\odot=\emptyset$.
Conversely, if $r_1$, $r_2$ and $r_3$ are subsets of $X\neq\emptyset$  such that (1)~$r_1\cup r_2\cup r_3=X$, (2)~$r_1\cap r_2=r_1 \cap r_3=r_2\cap r_3=\emptyset$, then there exists a unique $r\in {\bf 3}^X$  such that $r=( r_1, r_2, r_3)$.

\begin{remark} \label{ejP1LFI1} At this point, it is worth noting that, if we consider other 3-valued matrices, we can induce  in ${\bf 3}^X$ the corresponding operations. For example:
\begin{itemize}
\item[(1)] \ Let $\mm M_{\bf P1}=\langle\langle{\bf 3},\imp,\no\rangle,\{1,\frac{1}{2}\}\rangle$ be the 3-valued matrix associated to the Sette's logic $\bf P1$ (see~\cite{Sette}) where the truth-tables are dislayed below.
\begin{center}
$
\begin{array}{|c||c|c|c|}
\hline 
\imp & 1 & \frac 12 & 0 \\ 
\hline 
\hline 
1 & 1 & 1 & 0 \\ 
\hline 
\frac 12 & 1 & 1 & 0 \\ 
\hline 
0 & 1 & 1 & 1 \\ 
\hline 
\end{array} $ 
\quad
$
\begin{array}{|c|c|}
\hline 
 & \no \\ 
\hline 
1 & 0 \\ 
\hline 
\frac 12 & 1 \\ 
\hline 
0 & 1 \\ 
\hline 
\end{array} $
\end{center}
Then, for $r,u\in{\bf 3}^X$, we obtain the following:
\begin{itemize}
\item[] $\no r=(r_\ominus\cup r_\odot, r_\oplus,\emptyset)$
\item[]
$
r\imp u=
(r_\ominus\cup u_\oplus\cup u_\odot, \ (r_\oplus\cup r_\odot)\cap u_\ominus, \emptyset)
$
\end{itemize}
This coincides with the operations over triples for {\bf P1} proposed in~\cite[Definition~7.10.9]{CarCon}. 
\item[(2)] Let $\mm M_{\bf LFI1}=\langle\langle{\bf 3},\wedge,\vee,\to,\no,\cons\rangle,\{1,\frac{1}{2}\}\rangle$ be the 3-valued matrix associated to the propositional logic {\bf LFI1'} where the truth-tables are displayed below.
\begin{center}
$\begin{array}{|c||c|c|c|}
\hline 
\wedge & 1 & \frac 12 & 0 \\ 
\hline 
1 & 1 & \frac 12 & 0 \\ 
\hline 
\frac 12 & \frac 12 & \frac 12 & 0 \\ 
\hline 
0 & 0 & 0 & 0 \\ 
\hline 
\end{array} $
\quad
$ 
\begin{array}{|c||c|c|c|}
\hline 
\vee & 1 & \frac 12 & 0 \\ 
\hline 
1 & 1 & 1 & 1 \\ 
\hline 
\frac 12 & 1 & \frac 12 & \frac 12 \\ 
\hline 
0 & 1 & \frac 12 & 0 \\ 
\hline 
\end{array} $
\quad
$\begin{array}{|c||c|c|c|}
\hline 
\to & 1 & \frac 12 & 0 \\ 
\hline 
1 & 1 & \frac 12 & 0 \\ 
\hline 
\frac 12 & 1 & \frac 12 & 0 \\ 
\hline 
0 & 1 & 1 & 1 \\ 
\hline 
\end{array} $
\quad
$
\begin{array}{|c|c|c|}
\hline 
 & \no &\cons \\ 
\hline 
\hline 
1 & 0 &1\\ 
\hline 
\frac 12 & \frac 12 & 0\\ 
\hline 
0 & 1 & 1\\ 
\hline 
\end{array} $
\end{center}
As in the case of {\bf P1}, it can be seen that the operations over ${\bf 3}^X$ induced by the truth-tables above 
coincide with the operations over triples for {\bf LFI1'} proposed in~\cite[Definition~7.9.2]{CarCon} which, by its turn, were adapted from the operations over triples for {\bf MPT}  given in~\cite{consil}. As observed in~\cite{CarCon} and~\cite{consil}, the logics  {\bf LFI1'} and  {\bf MPT} coincide, up to language, with {\bf LFI1} and \dacdot. 
\end{itemize}
\end{remark}

\

\section{The paraconsistent first-order logic \quci} \label{s2}

In this section the logic \quci\ will be introduced, as a natural first-order version of \ciore, with semantics based on the notion of triples (or partial relations) analyzed in the previous section. Our approach is based on~\cite{car:etal:14} and~\cite[Chapter~7]{CarCon}. Hence, \quci\ will be defined as an axiomatic extension of \qmbc, the first-order version of \mbc\ proposed in~\cite{car:etal:14} and~\cite[Chapter~7]{CarCon}.

Recall that a first-order signature $\Sigma=\brac{\mm P,\mm F,\mm C}$ is composed of the following elements: a set $\mm P=\bigcup_{n\in\bb N}\mm P_n$ such that, for each $n\geq 1$, $\mm P_n$ is a set of predicate symbols of arity $n$; a set $\mm F=\bigcup_{n\in\bb N}\mm F_n$ such that, for each $n\geq 1$, $\mm F_n$ is a set of function symbols of arity $n$; and a set $\mm C$ of individual constants.
Let $\rr L(\Sigma)$ be the first-order language defined as  usual from the connectives $\wedge$, $\vee$, $\imp$, $\no$, $\cons$, the quantifiers $\pt$, $\ex$, a denumerable set $\mm V$ of individual variables and a given first-order signature $ \Sigma$. We denote by $Sent(\Sigma)$ the set of sentences (formulas without free-variables) over the signature $\Sigma$. If $\varphi$ is a formula, $x$ is a variable and $t$ is a term, then $\varphi(t/x)$ will denote the formula that results from $\varphi$ by replacing simultaneously all the free occurrences of the variable $x$ by $t$.

\begin{definition} Let $\Sigma$ be a first-order signature. The logic $\quci$ on the language $\rr L(\Sigma)$ is defined as the Hilbert calculus obtained by extending $\ciore$ (expressed in the language $\rr L(\Sigma)$) by adding the following:

\

\noindent{\bf Axiom schemata:}

\begin{myequation}
  \varphi(t/x)\imp\ex x\varphi, \, \mbox{ if $t$ is a term free for $x$ in } \varphi
   \end{myequation}
\vspace{-0.5cm}

\begin{myequation}
  \pt x\varphi\imp\varphi(t/x), \, \mbox{ if $t$ is a term free for $x$ in } \varphi
   \end{myequation}
\vspace{-0.5cm}

\begin{myequation}
\cons\ex x\varphi\imp\ex x\cons\varphi,   
\end{myequation}
\vspace{-0.5cm}

\begin{myequation}
\cons\pt x\varphi\imp\ex x\cons\varphi,  
   \end{myequation}
\vspace{-0.5cm}

\begin{myequation}
\exists x\cons\varphi\rightarrow\cons\exists x\varphi,   
   \end{myequation}
\vspace{-0.5cm}

\begin{myequation}
\exists x\cons\varphi\rightarrow\cons\pt x\varphi. 
   \end{myequation}
\vspace{-0.5cm}

\noindent{\bf Inference rules}
\begin{align*}
(\pt\mbox{-}{\bf In}) \ \frac{\alpha\imp\beta}{\alpha\imp\pt x\beta}
&& \mbox{if $x$ is not free in $\alpha$}\\
(\ex\mbox{-}{\bf In}) \ \frac{\alpha\imp\beta}{\ex x\alpha\imp\beta}
&& \mbox{if $x$ is not free in $\beta$}
\end{align*}
\end{definition}

\

\begin{theorem}\label{ciorequant} 
{\bf QCiore} proves the following schemas:

\begin{itemize}
\item[(i)] $\pt x\sneg\alpha\imp\sneg\ex x\alpha$,

\item[(ii)] $\sneg\ex x\alpha\imp \pt x\sneg \alpha$,

\item[(iii)] $\sneg\pt x\alpha\imp\ex x\sneg\alpha$,

\item[(iv)] $\big(\exists x\varphi\wedge\no\exists x\varphi\big)\rightarrow\forall x(\varphi\wedge\no\varphi)$,

\item[(v)] $\big(\forall x\varphi\wedge\no\forall x\varphi\big)\rightarrow\forall x(\varphi\wedge\no\varphi)$.
\end{itemize}
\end{theorem}
\begin{dem} 
Items (i)-(iii) can be proved from the similar results obtained for  {\bf QmbC} in~\cite[Chapter~7]{CarCon}. 

\noindent
(iv): 

\begin{tabular}{ll|l}
1. & $\exists x\cons\varphi\rightarrow\cons\exists x\varphi$ & {\bf (Ax16)} \\ 
2. & $\sneg\cons\exists x\varphi\rightarrow\sneg\exists x\cons\varphi$ & 1, Theorem \ref{esquema}(ii) \\ 
3. & $(\exists x\varphi\wedge\neg\exists x\varphi)\rightarrow\sneg\cons\exists x\varphi$ & Theorem \ref{esquema}(iii)\\ 
4. & $(\exists x\varphi\wedge\neg\exists x\varphi)\rightarrow\sneg\exists x\cons\varphi$ & 3, 2, transitivity of $\to$\\ 
5. & $\sneg\exists x\cons\varphi\rightarrow\forall x\sneg\cons\varphi$ & Item~(ii) \\ 
6. & $(\exists x\varphi\wedge\neg\exists x\varphi)\rightarrow\forall x\sneg\cons\varphi$ & 4, 5, transitivity of $\to$ \\ 
7. & $\forall x\sneg\cons\varphi\rightarrow\sneg\cons\varphi$ & {\bf (Ax12)} \\ 
8. & $\sneg\cons\varphi\rightarrow(\varphi\wedge\neg\varphi)$ 
& Theorem \ref{esquema}(iii)\\ 
9. & $\forall x\sneg\cons\varphi\rightarrow(\varphi\wedge\neg\varphi)$ & 7, 8, transitivity of $\to$ \\ 
10. & $\forall x\sneg\cons\varphi\rightarrow\forall x(\varphi\wedge\neg\varphi)$ & 9, {\bf ($\forall$-In)} \\ 
11. & $(\exists x\varphi\wedge\neg\exists x\varphi)\rightarrow\forall x(\varphi\wedge\neg\varphi)$ & 6, 10, transitivity of $\to$ 
\end{tabular} 

\

\noindent
(v): 

\begin{tabular}{ll|l}
1. & $\exists x\cons\varphi\rightarrow\cons\forall x\varphi$ & {\bf (Ax17)} \\ 
2. & $\sneg\cons\forall x\varphi\rightarrow\sneg\exists x\cons\varphi$ & 1, Theorem \ref{esquema}(ii) \\ 
3. & $(\forall x\varphi\wedge\neg\forall x\varphi)\rightarrow\sneg\cons\forall x\varphi$ &  Theorem \ref{esquema}(iii)\\ 
4. & $(\forall x\varphi\wedge\neg\forall x\varphi)\rightarrow\sneg\exists x\cons\varphi$ & 3, 2, transitivity of $\to$ \\ 
5. & $\sneg\exists x\cons\varphi\rightarrow\forall x\sneg\cons\varphi$ & Item~(ii) \\ 
6. & $(\forall x\varphi\wedge\neg\forall x\varphi)\rightarrow\forall x\sneg\cons\varphi$ & 4, 5, transitivity of $\to$ \\ 
7. & $\forall x\sneg\cons\varphi\rightarrow\sneg\cons\varphi$ & {\bf (Ax12)} \\ 
8. & $\sneg\cons\varphi\rightarrow(\varphi\wedge\neg\varphi)$ 
& Theorem \ref{esquema}(iii)\\ 
9. & $\forall x\sneg\cons\varphi\rightarrow(\varphi\wedge\neg\varphi)$ & 7, 8, transitivity of $\to$ \\ 
10. & $\forall x\sneg\cons\varphi\rightarrow\forall x(\varphi\wedge\neg\varphi)$ & 9, {\bf ($\forall$-In)} \\ 
11. & $(\forall x\varphi\wedge\neg\forall x\varphi)\rightarrow\forall x(\varphi\wedge\neg\varphi)$ & 6, 10, transitivity of $\to$  
\end{tabular} 

\end{dem}

\begin{theorem} \label{genrule}
The Generalization rule is derivable in {\bf QCiore}, i.e.
\begin{center}
{\bf (Gen)} \quad $\alpha \ \vdash_{\bf QCiore} \ \pt x \alpha$
\end{center}
\end{theorem}        
\begin{dem}

\begin{tabular}{ll|l}
1. & $\alpha$ & Hyp. \\ 
2. & $\alpha\imp\big((\pt x\alpha\imp\pt x\alpha)\imp\alpha\big)$ & {\bf (Ax 1)} \\ 
3. & $(\pt x\alpha\imp\pt x\alpha)\imp\alpha$ & 1,2, {\bf (MP)} \\ 
4.  & $(\pt x\alpha\imp\pt x\alpha)\imp\pt x\alpha$ & 3, {\bf ($\pt$-In)} \\ 
5.  & $\pt x\alpha\imp\pt x\alpha$ & Theorem \ref{esquema}(i) \\ 
6. & $\pt x\alpha$ & 4, 5, {\bf (MP)} \\ 
\end{tabular} 

\end{dem}

\

\
 
Finally, it is worth mentioning that it is possible to obtain a weak version of the {\em Deduction Meta-theorem} whose demonstration is identical to the one given for {\bf QmbC} in \cite[Chapter~7]{CarCon}.

\

\begin{theorem} \label{WDMT}(Weak  deduction Meta-Theorem for {\bf QCiore})
Suppose that there exists in {\bf QCiore} a derivation of $\psi$ from $\Gamma\cup\{\varphi\}$ such that no application of the rules {\bf ($\ex$-In)} and {\bf ($\pt$-In)} have, as their quantified variables, free variables of $\varphi$ (in particular, this holds when $\varphi$ is a sentence). Then $\Gamma\vdash_{\bf QCiore} \varphi\imp\psi$.
\end{theorem}

\

\section{Triples and twist structures for \ciore}

The so-called {\em twist structures} constitute an algebraic semantics for logics in terms of pairs of elements of a given class of algebras. The intuition behind twist structures  semantics is that a pair $(a,b)$ represents a truth-value for a formula $\varphi$, while $b$  represents a truth-value for the negation of $\varphi$. Twist structures  semantics were independently introduced by  M. Fidel~\cite{fid:78} and D. Vakarelov~\cite{vaka:77}, in order to obtain a representation of Nelson algebras in terms of twist structures over Heyting algebras. As observed by R.~Cignoli in~\cite{cin:86}, Fidel-Vakarelov structures are a particular case of a general construction proposed by J. Kalman in 1958 (see~\cite{kal:58}). 
In~\cite[Definition~9.18]{ConFigGol} M. Coniglio, A. Figallo-Orellano and A. Golzio presented the following twist structures semantics for \ciore:

\begin{definition}  \label{defKCiore}
Let $\mathcal{A} = \langle A, \sim , \supset, \sqcup, \sqcap, 0, 1\rangle$ be a Boolean algebra and let $\mathbb{P}_{\mathcal{A}} = \{(a, b) \in A \times A \ : \ a \sqcup b = 1\}$.
A {\em twist structure for \ciore\ over $\mathcal{A}$}  is an algebra $\mathcal{P}_\mathcal{A}$ over the  signature $\{\wedge,\vee,\imp,\no,\cons\}$ with domain $\mathbb{P}_{\mathcal{A}}$ such that the operations are defined as follows, for every $(z_1,z_2),(w_1,w_2) \in \mathbb{P}_{\mathcal{A}}$:

\begin{itemize}
 \item[(i)] $(z_1,z_2)\wedge (w_1,w_2) = (z_1\sqcap w_1,(z_1 \sqcap w_1) \supset ((z_1 \sqcap z_2) \sqcap (w_1 \sqcap w_2)))$;
 \item[(ii)] $(z_1,z_2)\vee (w_1,w_2) = (z_1\sqcup w_1,(z_1 \sqcup w_1) \supset ((z_1 \sqcap z_2) \sqcap (w_1 \sqcap w_2)))$;
 \item[(iii)] $(z_1,z_2) \to (w_1,w_2) = (z_1 \supset w_1,(z_1 \supset w_1) \supset ((z_1 \sqcap z_2) \sqcap (w_1 \sqcap w_2)))$;
\item[(iv)] $\neg (z_1,z_2) = (z_2,z_1)$;
 \item[(v)] $\circ (z_1,z_2) = ({\sim}(z_1 \sqcap z_2),z_1 \sqcap z_2)$.\\[-2mm]
\end{itemize}
\end{definition}

\noindent
In what follows, it will be proven that the twist structures for \ciore, which are algebras of pairs, are equivalent to the algebras of triples presented in Section~\ref{s1.5}. In order to do this, it will be observed that each algebra ${\bf 3}^X$, which is formed by triples over the Boolean algebra $\wp(X)$ (by Proposition~\ref{lfi1al}), can be generalized to any Boolean algebra.

\begin{definition} \label{defTCiore}
Let $\mathcal{A}$ be a Boolean algebra and let $\mathbb{T}_{\mathcal{A}} = \{(a, b,c) \in A \times A \times A \ : \ a \sqcap b =a \sqcap c =b \sqcap c =0$, and $a \sqcup b \sqcup c= 1\}$.
A {\em  structure of triples for \ciore\ over $\mathcal{A}$}  is an algebra $\mathcal{T}_\mathcal{A}$ over the  signature $\{\wedge,\vee,\imp,\no,\cons\}$ with domain $\mathbb{T}_{\mathcal{A}}$ such that the operations are defined as in  Proposition~\ref{lfi1al}, by changing $\cap$ and $\cup$ by $\sqcap$ and $\sqcup$, respectively.
\end{definition}

\noindent Thus, if  $z=(z_1,z_2,z_3)$ and $w=(w_1,w_2,w_3)$ are in $\mathbb{T}_{\mathcal{A}}$ then, for instance, $\no z = (z_2,z_1,z_3)$ and $z\imp w = ( z_2\sqcup w_1,  (z_1\sqcup z_3)\sqcap w_2, (z_1\sqcup z_3)\sqcap w_3)$. It is clear that, given a non-empty set $X$,  the algebra $\mathbb{T}_{\wp(X)}$ is isomorphic to  ${\bf 3}^X$, by  Proposition~\ref{lfi1al}. Hence, Definition~\ref{defTCiore} generalizes the algebras of triples  ${\bf 3}^X$ from $\wp(X)$  to any Boolean algebra $\mathcal{A}$.

\begin{remark} \label{obsbola} Let $\mathcal{A}$ be a Boolean algebra.
For every $z \in \mathbb{P}_{\mathcal{A}}$ and every  $w \in \mathbb{T}_{\mathcal{A}}$ let $\bot^P_\mathcal{A} = (z \wedge \neg z) \wedge \circ z$ and $\bot^T_\mathcal{A} = (w \wedge \neg w) \wedge \circ w$. By Definitions~\ref{defKCiore} and~\ref{defTCiore} it is easy to see that $\bot^P_\mathcal{A} =(0,1)$ and $\bot^T_\mathcal{A} =(0,1,0)$, hence the choice of $z$ and $w$ for defining $\bot^P_\mathcal{A}$ and $\bot^T_\mathcal{A}$ is irrelevant. Moreover, $\circ z = (z \wedge \neg z)  \to \bot^P_\mathcal{A}$ and  $\circ w= (w \wedge \neg w)  \to \bot^T_\mathcal{A}$ for every $z \in \mathbb{P}_{\mathcal{A}}$ and every $w \in \mathbb{T}_{\mathcal{A}}$. 
\end{remark}

\begin{proposition} \label{iso-par-tri}
Let $\mathcal{A}$ be a Boolean algebra. Then, the function $\dagger:\mathbb{T}_{\mathcal{A}} \to\mathbb{P}_{\mathcal{A}}$ given by $\dagger(z_1,z_2,z_3)=(z_1 \sqcup z_3,z_2 \sqcup z_3)$ is an isomorphism between $\mathcal{T}_\mathcal{A}$ and $\mathcal{P}_\mathcal{A}$. The inverse of $\dagger$ is given by $\ddagger:\mathbb{P}_{\mathcal{A}} \to\mathbb{T}_{\mathcal{A}}$ such that $\ddagger(z_1,z_2)=(z_1\sqcap {\sim}z_2,z_2\sqcap {\sim}z_1,z_1 \sqcap z_2)$.  
\end{proposition}
\begin{dem}
Straightforward. Observe that, by Remark~\ref{obsbola}, the fact that $\dagger(\circ z)={\circ}{\dagger}(z)$ follows from the fact that $\dagger(\neg z)=\neg{\dagger}(z)$, $\dagger(z \div w)={\dagger}(z) \div {\dagger}(w)$ (for $\div \in\{\wedge,\to\}$) and ${\dagger}(\bot^T_\mathcal{A})=\bot^P_\mathcal{A}$. The details of the proof are left to the reader.
\end{dem}

\

\noindent The last result shows that the algebras of triples ${\bf 3}^X$ defined in Section~\ref{s1.5} correspond to twist structures for \ciore\ over Boolean algebras of the form $\wp(X)$, as defined in~\cite{ConFigGol}.

\

\section{\quci-structures, satisfiability, validity.}\label{s3}

Inspired by Proposition \ref{lfi1al}, the notion of \quci-structure will be introduced in this section. As it will be shown below, these  first-order structures will be adequate for characterizing the logic \quci.

A {\em partial (or pragmatic) structure over the signature $ \Sigma$}  is an ordered pair $\ff A=\langle A,(\cdot)^\ff A\rangle$ 
where $A\neq\emptyset$ and $(\cdot)^\ff A$ is a function such that

\begin{itemize}
\item[-] for all $R\in\mm P_n$, 
$R^{\ff A}=\langle R^{\ff A}_{ \oplus},R^{\ff A}_{ \ominus},R^{\ff A}_{ \odot}\rangle$
is a triple  over $A^n$,
\item[-] $(\cdot)^\ff A$ is defined as usual over $\mm F$ and $\mm C$.
\end{itemize}

If 
$\Sigma=\langle\mm P,\mm F,\mm C\rangle$ and 
$\Sigma'=\langle\mm P',\mm F',\mm C'\rangle$ 
are first-order signatures such that 
$\mm P'_n\subseteq\mm P_n$,
$\mm F'_n\subseteq\mm F_n$ for every $n \geq 1$, and
$\mm C'\subseteq\mm C$, we write $\Sigma'\subseteq\Sigma$. 
In this case if $\ff A=\langle A,(\cdot)^{\ff A}\rangle$ is a \quci-structure for $\Sigma$ and $\ff A'=\langle A,(\cdot)^{\ff A'}\rangle$ is a \quci-structure for $\Sigma'$ such that $(\cdot)^{\ff A'}$ is the restriction of the map $(\cdot)^\ff A$ over $\Sigma'$ we say that $\ff A'$ is the {\em reduct} of $\ff A$ to $\Sigma'$, denoted by $\ff A|\Sigma'=\ff A'$.

Next, we recall some notions that are standard for Tarskian structures, as well as the notion of {\em quantification over sets}  considered by  Coniglio and Silvestrini in \cite{consil}. 

\begin{definition}\label{def-assi}
Let $A$ be a non-empty set. An {\em assignment into $A$} is any map $s:\mm V\rightarrow  A$. We denote by $S(A)$ the set of all assignments into $A$, i.e. $S(A)=A^\mm V$. If $\ff A$ is a partial structure over $\Sigma$ with domain $A$ then an {\em assignment into $\ff A$} is any assignment into $A$. The set of all assignments into $\ff A$ will be denoted by $S(\ff A)$, i.e. $S(\ff A)=S(A)=A^\mm V$.
\end{definition}

\begin{definition}\label{def-valterm}
Let $s \in S(\ff A)$. The {\em value of the term $t$ in $\ff A$ under the assignment $s$}, denoted by $t^\ff A[s]$, is defined inductively as follows:

\begin{itemize}
\item[-] if $t$ is $v$, for $v\in\mm V$, then $t^\ff A[s]:=s(v)$,
\item[-] if $t$ is $c$, for $c\in\mm C$, then $t^\ff A[s]:=c^\ff A$,
\item[-] if $t$ is $f(t_1,\ldots,t_n)$, where $f\in\mm F_n$ and $t_1,\ldots,t_n$  are terms, then $t^\ff A[s]:=f^\ff A(\ter{1},\ldots,\ter{n})$.
\end{itemize}
\end{definition}

\begin{definition}\label{def-sxa}
Let $A$ be a non-empty set. Given $s \in S(A)$ $x\in\mm V$ and $a\in A$, the assignment $s_x^a:\mm V\rightarrow A$ is defined by

\begin{equation}\tag{1}\label{eq01}
s_x^a(y)=\left\{
\begin{array}{ll}
s(y)&\mbox{if}\ y\neq x\\
a&\mbox{if}\ y=x
\end{array}
\right.
\end{equation}
\end{definition}

\begin{definition}\label{def-quant}
Let $A$ be a non-empty set. Given a set $Z$ let $\wp(Z)$  and $\wp(Z)_+$ be the set of subsets and of  non-empty subsets of $Z$, respectively. For $x\in\mm V$ let $\widehat{\pt x}:\wp(S(A)) \to  \wp(S(A))$ and $\widehat{\ex x}:\wp(S(A)) \to  \wp(S(A))$ be defined as follows, for every $Y \subseteq S(A)$:

$${\widehat{\pt x}(Y)}:=\{s\in S(A) \ : \  s_x^a\in Y\mbox{ for all }a\in A\},$$
$${\widehat{\ex x}(Y)}:=\{s\in S(A) \ : \  s_x^a\in Y\mbox{ for some }a\in A\}.$$
The functions $\widetilde{\forall}:\wp({\bf 3})_+ \to {\bf 3}$ and  $\widetilde{\exists}:\wp({\bf 3})_+ \to {\bf 3}$ are defined as follows, for every $\emptyset\neq Y \subseteq {\bf 3}$:
$$
\widetilde{\forall}(Y):=
\left\{
\begin{array}{ll}
1
& \ \mbox{if }\ 
1 \in Y, \ 0 \notin Y
\\[1mm]
\frac{1}{2}
& \ \mbox{if }\ 
Y=\{\frac{1}{2}\}
\\[1mm]
0
& \ \mbox{if }\ 
0 \in Y
\end{array}
\right.
\hspace{2cm}
\widetilde{\exists}(Y):=
\left\{
\begin{array}{ll}
1
& \ \mbox{if }\ 
Y\neq\{\frac{1}{2}\}, \ Y\neq\{0\}
\\[1mm]
\frac{1}{2}
& \ \mbox{if }\ 
Y=\{\frac{1}{2}\}
\\[1mm]
0
& \ \mbox{if }\ 
Y=\{0\}
\end{array}
\right.
$$
\end{definition}

\

\begin{definition} \label{defLFI2structure} A {\em \quci-structure} over the signature $\Sigma$ is a pair
$\langle\ff A,\tripc{\cdot}\rangle$ such that $\ff A$ is a partial structure and  $\tripc{\cdot}:{\rm For}(\Sigma)\imp {\bf 3}^{S(\ff A)}$ is a map defined recursively as follows, for every $s \in S(\ff A)$ (recall the operations defined in Remark~\ref{op3X} and Definition~\ref{def-quant}):
 \begin{enumerate}
 \item If $P(t_1,\ldots,t_n)$ is an atom, $\tripc{P(t_1,\ldots,t_n)}(s)=
P^\ff A(t_1^{\ff A}[s],\ldots,t_n^{\ff A}[s])$,

\item $\tripc{\varphi\wedge\psi}=\tripc{\varphi}\wedge \tripc{\psi}$,

\item $\tripc{\varphi\vee\psi}=\tripc{\varphi}\vee \tripc{\psi}$,

\item $\tripc{\varphi\imp\psi}=\tripc{\varphi} \imp \tripc{\psi}$,

\item 
$\tripc{\no \varphi}=\no \tripc{\varphi}$,

\item 
$\tripc{\cons\varphi}=\cons\tripc{\varphi}$,

\item 
$\tripc{\pt x \varphi}(s)=\widetilde{\forall}(\{\tripc{\varphi}(s_x^a) \ : \ a\in A\})$. 

\item
$\tripc{\ex x \varphi}(s)=\widetilde{\exists}(\{\tripc{\varphi}(s_x^a) \ : \ a\in A\})$.
 
 \end{enumerate}
\end{definition}

\

\begin{proposition} \label{charLFI2structure} Let 
$\langle\ff A,\|\cdot\|^\ff A\rangle$ be a {\em \quci-structure} over the signature $\Sigma$. Then, by denoting $\tripc{\varphi}$ as $\langle \mas\varphi, \menos\varphi, \por\varphi \rangle$ (as stipulated in Section~\ref{s1.5}), the following holds (recalling Remark~\ref{op3X}):
 \begin{itemize}
 
 \item[-] If $P(t_1,\ldots,t_n)$ is an atom, then, for every $\# \in \{\oplus,\ominus,\odot\}$ and $s \in S(\ff A)$:
\begin{center}
\begin{tabular}{rcl}
$s\in \|P(t_1, \ldots,t_n)\|^{\mathfrak A}_\#$
& \ \mbox{ iif }\ 
&$(t_1^{\ff A}[s],\ldots,t_n^{\ff A}[s])\in P^{\mathfrak A}_\#$,
\end{tabular}
\end{center} 
 
\item[-] 
$\tripc{\varphi\wedge\psi} \ := \ 
\left\langle
(\mas\varphi\cup\mas\psi)\cap
(\mas\varphi\cup\por\psi)\cap
(\por\varphi\cup\mas\psi),\right.$

$\phantom{\tripc{\varphi\wedge\psi} \ := \ } \ \
\left.\menos\varphi\cup\menos\psi,
\por\varphi\cap\por\psi
\right\rangle
$
,

\item[-] 
$\tripc{\varphi\vee\psi} \ := \ 
\left\langle
\mas\varphi\cup\mas\psi\cup
(\menos\varphi\cap\por\psi)\cup
(\por\varphi\cap\menos\psi), \right.$

$\phantom{\tripc{\varphi\vee\psi} \ := \ } \ 
\left.\menos\varphi\cap\menos\psi,
\por\varphi\cap\por\psi
\right\rangle
$,

\item[-] 
$\tripc{\varphi\imp\psi} \ := \ \left\langle
\menos\varphi\cup
(\mas\varphi\cap\mas\psi)\cup
(\mas\varphi\cap\por\psi)\cup
\right.$

$\phantom{\tripc{\varphi\imp\psi}\ := \ }
\left.
(\por\varphi\cap\mas\psi),
(\mas\varphi\cup\por\varphi)\cap\menos\psi,
\por\varphi\cap\por\psi
\right\rangle
$,

\item[-] 
$\tripc{\no \varphi}\ := \ 
\langle
\menos\varphi,
\mas\varphi,
\por\varphi
\rangle$,

\item[-]
$\tripc{\cons\varphi}\ :=\ 
\langle
\mas\varphi\cup\menos\varphi,
\por\varphi,
\emptyset
\rangle$,

\item[-] 
$\tripc{\pt x \varphi}\ := 
\left\langle 
\widehat{\ex x}(\mas{\varphi})-\widehat{\ex x}(\menos\varphi),
\widehat{\ex x}(\menos{\varphi}),
\widehat{\pt x}(\por{\varphi})
\right\rangle$,
\
\item[-] 
$\tripc{\ex x \varphi}\ := 
\left\langle 
S(\ff A)-\big(\widehat{\pt x}(\menos\varphi)\cup\widehat{\pt x}(\por\varphi)\big),\ 
\widehat{\pt x}(\menos{\varphi}),\ 
\widehat{\pt x}(\por{\varphi})
\right\rangle$.
 \end{itemize}
\end{proposition}
\begin{dem}
Straightforward from the definitions and from Proposition~\ref{lfi1al}.
\end{dem}

\

\begin{remark} \label{obser}
Observe that, with notation as in Proposition~\ref{charLFI2structure}: $\mas{\varphi} = \{ s \in S(\ff A)   \ : \  \tripc{\varphi}(s)=1\}$;  $\por{\varphi} = \{ s \in S(\ff A)   \ : \  \tripc{\varphi}(s)=\frac{1}{2}\}$; and  $\menos{\varphi} = \{ s \in S(\ff A)  \ : \  \tripc{\varphi}(s)=0\}$, for every formula $\varphi$.
\end{remark}

\

\begin{remark}
Recall the structure of triples for \ciore\  introduced in Definition~\ref{defTCiore}. By considering the Boolean algebra $\wp(S(\ff A))$, Proposition~\ref{charLFI2structure} states operations $[\forall x]_T^{\ff A},[\exists x]_T^{\ff A}:\mathcal{T}_{\wp(S(\ff A))} \to \mathcal{T}_{\wp(S(\ff A))}$ (for $x \in \mm V$) given by 
$$[\forall x]_T^{\ff A}\langle Z_1,Z_2,Z_3\rangle=\left\langle 
\widehat{\ex x}(Z_1)-\widehat{\ex x}(Z_2), \
\widehat{\ex x}(Z_2), \
\widehat{\pt x}(Z_3)
\right\rangle,$$
$$[\exists x]_T^{\ff A}\langle Z_1,Z_2,Z_3\rangle=\left\langle 
S(\ff A)-\big(\widehat{\pt x}(Z_2)\cup\widehat{\pt x}(Z_3)\big),\ 
\widehat{\pt x}(Z_2),\ 
\widehat{\pt x}(Z_3)
\right\rangle.$$ 
By using the isomorphism $\dagger$ between  $\mathcal{T}_\mathcal{A}$ and $\mathcal{P}_\mathcal{A}$ given in Proposition~\ref{iso-par-tri}, as well as its inverse $\ddagger$, the operations $[\forall x]_T^{\ff A}$ and $[\exists x]_T^{\ff A}$ induce operations $[\forall x]_P^{\ff A},[\exists x]_P^{\ff A}:\mathcal{P}_{\wp(S(\ff A))} \to \mathcal{P}_{\wp(S(\ff A))}$ as follows: $[\forall x]_P^{\ff A}\langle Z_1,Z_2\rangle:={\dagger}([\forall x]_T^{\ff A}({\ddagger}\langle Z_1,Z_2\rangle))$ and 
$[\exists x]_P^{\ff A}\langle Z_1,Z_2\rangle:={\dagger}([\exists x]_T^{\ff A}({\ddagger}\langle Z_1,Z_2\rangle))$. Observe that
$$[\forall x]_P^{\ff A}\langle Z_1,Z_2\rangle=\left\langle 
\big(\widehat{\ex x}(Z_1-Z_2)-\widehat{\ex x}(Z_2-Z_1)\big)\cup \widehat{\pt x}(Z_1 \cap Z_2),\
\widehat{\ex x}(Z_1 -Z_2) \cup \widehat{\pt x}(Z_1 \cap Z_2)
\right\rangle,$$
$$[\exists x]_P^{\ff A}\langle Z_1,Z_2\rangle=\left\langle 
\big(S(\ff A)-\widehat{\pt x}(Z_2-Z_1)\big)\cup\widehat{\pt x}(Z_1 \cap Z_2),\ 
\widehat{\pt x}(Z_2-Z_1) \cup 
\widehat{\pt x}(Z_1 \cap Z_2)
\right\rangle.$$
This shows how to define quantification over twist structures for \ciore, allowing to define twist structures semantics for \quci. Observe that, by definition, ${\dagger}([\forall x]_T^{\ff A}\langle Z_1,Z_2,Z_3\rangle) = [\forall x]_P^{\ff A}({\dagger}\langle Z_1,Z_2,Z_3\rangle)$ and ${\dagger}([\exists x]_T^{\ff A}\langle Z_1,Z_2,Z_3\rangle) = [\exists x]_P^{\ff A}({\dagger}\langle Z_1,Z_2,Z_3\rangle)$. That is, the isomorphism $\dagger$ also preserves quantification.
\end{remark}

\

\begin{definition} \label{defsat}
Let $\langle\ff A,\tripc{\cdot}\rangle$ be a \quci-structure and let $\varphi$ be a formula.  
We say that $s\in\ese A$ {\em \quci-satisfies} $\varphi$ in $\ff A$, denoted by $\ff A\dodash\varphi[s]$, if $s\in\mas\varphi\cup\por\varphi$. Besides,  
$\varphi$ is said to be {\em \quci-valid} in $\ff A$, denoted by $\ff A\Vdash_\quci\varphi$, if  $\mas{\varphi}\cup\por{\varphi}\ =\ S(\ff A)$. A \quci-structure $\ff A$ is a \quci-model of a set $\Gamma$ of formulas if $\ff A\dodash\varphi$ for each $\varphi\in\Gamma$.

Given a set $\Gamma\cup \{ \varphi \}$ of formulas, we say that {\em $\varphi$ is a \quci-consequence of $\Gamma$}, denoted by $\Gamma\Vdash_\quci\varphi$ if, for every
\quci-structure $\ff A$, we have that: $\ff A\Vdash_\quci\psi$ for every $\psi\in \Gamma$ implies that $\ff A\Vdash_\quci\varphi$.
A formula $\varphi$ is {\em \quci-valid}, denoted by $\Vdash_\quci \varphi$, if  $\emptyset\Vdash_\quci\varphi$.
\end{definition}

\

\noindent The following propositions will be useful in the sequel.

\

\begin{proposition}\label{s3.6c} Let $\langle\ff A,\tripc{\cdot}\rangle$ be a \quci-structure and let $\varphi$ and $\psi$ be formulas. Then
 \begin{itemize}
 \item[(i)] $\ff A \dodash P(t_1, \ldots,t_n)[s]$ \ iff \ $(\ter{1},\ldots,\ter{n})\in\tas{P}\cup\tpor{P}$, \, for every atomic formula  $P(t_1, \ldots,t_n)$.
 \item[(ii)] $\ff A\dodash(\varphi\wedge\psi)[s]$ \ iff \ $\ff A\dodash\varphi[s]$ and $\ff A\dodash\psi[s]$.
 \item[(iii)] $\ff A\dodash(\varphi\vee\psi)[s]$ \ iff \ $\ff A\dodash\varphi[s]$ or $\ff A\dodash\psi[s]$.
 \item[(iv)] $\ff A\dodash(\varphi\imp\psi)[s]$ \ iff \ $\ff A\dodash\varphi[s]$ implies $\ff A\dodash\psi[s]$.
 \item[(v)] $\ff A\dodash\no\varphi[s]$ \ iff \ $s\in\menos{\varphi}\cup\por{\varphi}$.
 \item[(vi)] $\ff A\dodash\cons\varphi[s]$ \ iff \ $s\in\mas{\varphi}\cup\menos{\varphi}$.
 \item[(vii)]  $\ff A\dodash\pt x \varphi[s]$ \ iff \ $\ff A\dodash \varphi[s_x^a]$ for all $a \in A$.
 \item[(viii)] $\ff A\dodash\ex x \varphi[s]$ \ iff \ there exists $a \in A$ such that $\ff A\dodash \varphi[s_x^a]$.
 \end{itemize}
\end{proposition}
\begin{dem} It is consequence of Definition \ref{defLFI2structure}.
\end{dem}

\

\begin{proposition}\label{cis3.8}
Let $\varphi$ a formula, $\ff A$ a \quci-structure and $s\in S(\ff A)$. Then:
\begin{itemize}
\item[(i)] $s\in\mas{\varphi}$ \ iff \ $\ff A\dodash(\varphi\wedge\cons\varphi)[s]$.
\item[(ii)] $s\in\menos{\varphi}$ \ iff \ $\ff A\dodash(\no\varphi\wedge\cons\varphi)[s]$.
\item[(iii)] $s\in\por{\varphi}$ \ iff \ $\ff A\dodash(\varphi\wedge\no\varphi)[s]$.
\item[(iv)] $s$ \quci-satisfies in $\ff A$ at most two of the formulas \ $\varphi, \ \no\varphi,\ \cons\varphi.$
\end{itemize}
\end{proposition}
\begin{dem} (i) $\ff A\dodash(\varphi\wedge\cons\varphi)[s]$ iff, by Proposition~\ref{s3.6c}~(ii), 
$\ff A\dodash\varphi[s]$ and $\ff A\dodash\cons\varphi[s]$ iff, by Definition~\ref{defsat} and  Proposition~\ref{s3.6c}~(vi), $s\in\mas{\varphi}\cup\por{\varphi}$ and $s\in\mas{\varphi}\cup\menos{\varphi}$ iff $s\in\mas{\varphi}$, since $\por{\varphi} \cap \menos{\varphi}=\emptyset$.\\
Items (ii), (iii) and (iv) are proved analogously.
\end{dem}

\

\section{Soundness and Completeness}\label{s4}

In this section, we shall prove the soundness and completeness of \quci\ with respect to its first-order structures. For this, it is necessary to establish some auxiliary results which we shall use repeatedly in what follows. From now on $\langle\ff A,\tripc\cdot\rangle$ denotes a \quci-structure.

\begin{proposition}\label{citad.imp} 
The following statements are equivalent:
\begin{itemize}
\item[(i)] $\ff A\dodash\varphi\imp\psi$;
\item[(ii)] $\mas{\varphi}\cup\por{\varphi}\subseteq \mas{\psi}\cup\por{\psi}$;
\item[(iii)] $\menos{\psi}\subseteq \menos{\varphi}$.
\end{itemize}
\end{proposition}
\begin{dem} It is consequence of Definition~\ref{defLFI2structure} and Proposition~\ref{s3.6c}~(iv).
\end{dem}

\

\noindent Next we shall prove that the inference rules of {\bf Qciore} preserve validity.

\begin{proposition}\label{cim.tres}
If $\ff A\dodash\varphi$ and $\ff A\dodash\varphi\imp\psi$, then $\ff A\dodash\psi$.
\end{proposition}
\begin{dem} Since $\mas\varphi \cup \por\varphi=S(\ff A)$ and taking into account that $\mas\varphi\cup\por\varphi\subseteq\mas\psi\cup\por\psi$, it follows that $\ff A\dodash\psi$.
\end{dem}

\

An {\em instance} of a propositional formula $\varphi=\varphi(p_1,\ldots,p_n)$ in $\rr L$ is a formula in $\rr L(\Sigma)$ obtained from $\varphi$ by simultaneously substituting each occurrence of the propositional variable $p_i$ by the formula $\beta_i$ in $\rr L(\Sigma)$, where $1 \leq i \leq n$. We denote such instance by $\varphi[\beta_1/p_1,\ldots,\beta_n/p_n]$.

\

\begin{proposition}\label{cim.siete}
Let $\alpha=\varphi[\beta_1/p_1,\ldots,\beta_n/p_n]$ be an instance of a propositional formula $\varphi(p_1,\ldots,p_n)$ in $\rr L_\Theta$. For every $s\in S(\ff A)$ we define the valuation $v_s$ over the 3-valued matrix $M_e$ for \ciore\ (recall Theorem~\ref{SoundCompFI2}) as follows:

$$
v_s(p_i)=\left\{
\begin{array}{ll}
1     &\ \mbox{ if }     \  s\in \mas{\beta_i}    \\[2mm]
\frac{1}{2}     &\ \mbox{ if }     \  s\in \por{\beta_i}    \\[2mm]
0     &\ \mbox{ if }     \  s\in \menos{\beta_i}    
\end{array}
\right.
$$

\noindent for $1 \leq i\leq n$, and $v_s(p_i)$ is arbitrary otherwise. Then $\tripc{\alpha}(s)=v_s(\varphi)$.
\end{proposition}
\begin{dem} Straightforward by induction on the complexity of $\varphi$, taking into account Definition~\ref{defLFI2structure} and Remark~\ref{obser}.
\end{dem}

\

\noindent As a direct consequence, we have the following

\begin{corollary}\label{cim.sietee}
If $\varphi$ is a instance of  a \eledos -tautology, then $\ff A\dodash\varphi$.							
\end{corollary}

\begin{proposition}\label{ci1.15}
Let $t(x_1,\ldots,x_n)$ be a term  and let $\varphi(x_1,\ldots,x_n)$ be a formula in $\rr L(\Sigma)$.
If $s$ and $s'$ are assignments in $\ff A$ such that $s(x_i)=s'(x_i)$ for all $i$ ($1\leq i\leq n$), then
\begin{itemize}
\item[(i)] $t^\ff A[s]=t^\ff A[s']$,
\item[(ii)] $\tripc\varphi(s)=\tripc\varphi(s')$,
\item[(iii)] $\ff A\dodash\varphi[s]$ \ iff \ $\ff A\dodash\varphi[s']$.
\end{itemize}
\end{proposition}
\begin{dem} Straightforward, by induction on the length of the term $t$ and the complexity of  $\varphi$.
\end{dem}

\

\begin{theorem}\label{cim.nueve}
Let  $\varphi$ be a sentence. Then: either  $\mas{\varphi}=S(\ff A)$ or $\menos{\varphi}=S(\ff A)$ or  $\por{\varphi}=S(\ff A)$.
\end{theorem}
\begin{dem} Let $\varphi$ be a sentence and let $s\in S(\ff A)$. Suppose that $s\in\mas{\varphi}$. Then, for all  $s'\in S(\ff A)$  it holds  vacuously that $s(x)=s'(x)$ for every variable $x$ occurring free in $\varphi$. Then $s'\in\mas{\varphi}$, by Proposition~\ref{ci1.15}, and so $\mas{\varphi}=S(\ff A)$. Analogously,  $s\in\menos{\varphi}$ implies that $\menos{\varphi}=S(\ff A)$, and $s\in\por{\varphi}$ implies that $\por{\varphi}=S(\ff A)$.
\end{dem}

\

\begin{remark} Observe that if $\varphi$ is a sentence then there exists $\bf v\in \bf 3$ such that $\tripc\varphi(s)=\bf v$ for all $s\in\ese A$. In this case, we shall write $\tripc\varphi=\bf v$.
\end{remark}

\begin{corollary} \label{corosent}
If $\varphi$ is a sentence then one and only one of the following statements holds: \\[1mm]
(i) $\ff A\Vdash_\quci\varphi\wedge\cons\varphi$,\\
(ii)  $\ff A\Vdash_\quci\no\varphi\wedge\cons\varphi$,\\
(iii)  $\ff A\Vdash_\quci\varphi\wedge\no\varphi$.
\end{corollary}

\

\noindent Let $t$ and $u$ be terms and $s\in S(\ff A)$. Let $t'$ be the term that results from $t$ by replacing simultaneously all the occurrences of the variable $x$ by $u$, which will be denoted by $t(u/x)$. Let $s_{x}^{a}$ the assignment defined as in Definition~\ref{def-sxa}, where $a=u^\ff A[s]$ (recall Definition~\ref{def-valterm}). It is easy to prove, by induction on the complexity of the term $t$, that \, $(t')^\ff A[s] = t^\ff A [s_{x}^{a} ]$. Moreover:

\

\begin{proposition}\label{cim.lemasx} 
Let $t$ be a term free for $x$ in $\varphi$, and let $a=t^\ff A[s]$. Then:
\begin{itemize}
\item[(i)] $\tripc{\varphi(t/x) }(s) = \tripc\varphi \left(s_x^{a}\right)$ \ \ {\em (Substitution Lemma)},
\item[(ii)] $\ff A\dodash\varphi(t/x)[s]$ \ iff \ $\ff A\dodash\varphi\left[ s_x^{a}\right]$,
\item[(iii)] If $\ff A\dodash\pt x \varphi[s]$, \, then $\ff A\dodash\varphi(t/x)[s]$.
\end{itemize}
\end{proposition}
\begin{dem} It is routine.
\end{dem}

\

\begin{proposition}\label{cim.diez}
\begin{itemize}
\item[]
\item[(i)] Let $t$ be a term free for $x$ in $\varphi$. Then: \\[2mm]
\begin{tabular}{rlrl}
(a) & $\dodash\pt x \varphi\imp\varphi(t/x)$,  &  (b)  & $\dodash\varphi(t/x)\imp\ex x\varphi$.\\
\end{tabular}
\item[(ii)]  $\dodash\ex x\cons\varphi\sii\cons\pt x\varphi$, 
\item[(iii)] $\dodash\ex x\cons\varphi\sii\cons\ex x\varphi$.	
\end{itemize}
\end{proposition}
\begin{dem} (i) They are consequence of Lemma \ref{cim.lemasx}.\\[1mm]
(ii) It follows from $\menos{\cons\pt x\varphi}=\por{\pt x\varphi}=\widehat{\pt x}(\por\varphi)=\widehat{\pt x}(\menos{\cons\varphi})=\menos{\ex x\cons\varphi}$. \\[1mm]
(iii) It follows from $\menos{\ex x\cons\varphi}=\widehat{\pt x}(\menos{\cons\varphi})=\widehat{\pt x}(\por\varphi)=\por{\ex x\varphi}=\menos{\cons\ex x\varphi}$.
\end{dem}

\
 
\begin{proposition}\label{cim.once} 
\begin{itemize}
\item[]
\item[(i)] \ $\varphi\imp\psi\dodash\varphi\imp\pt x \psi$ \
if the variable $x$ does not occur free in $\varphi$. 

\item[(ii)] \ $\varphi\imp\psi\dodash\ex x\varphi\imp\psi$ \
if the variable $x$ does not occur free in $\psi$. 
\end{itemize}
\end{proposition}
\begin{dem} Let $\ff A$ be a \quci-structure and suppose that $\dodash\varphi\imp\psi$. Then, by Proposition \ref{citad.imp},  $\menos\psi\subseteq\menos\varphi$. \\[1mm]
 (i) Since $x$ does not occur free in $\varphi$ and from Proposition \ref{ci1.15}, it follows that $\widehat{\ex x}(\menos{\varphi})=\menos{\varphi}$. Since $\widehat{\ex x}(\menos\psi)\subseteq\widehat{\ex x}(\menos\varphi)$, it follows that $\menos{\pt x\psi}\subseteq\menos\varphi$.\\[1mm]
(ii) Since $x$ does not occur free in $\psi$ and from Proposition \ref{ci1.15}, it follows that $\widehat{\pt x}(\menos{\psi})=\menos{\psi}$. 
Since $\widehat{\pt x}(\menos\psi)\subseteq\widehat{\pt x}(\menos\varphi)$, it follows that $\menos{\psi}\subseteq\menos{\ex x\varphi}$.
\end{dem}

\

\noindent
At this point, it is important to remark that, unlike to what happens in the first-order version of {\bf LFI1}, the quantifiers $\pt$ and $\ex$ are not interdefinable in \quci\ in terms of the negation $\neg$. To see this, consider the following:  

\begin{remark} \label{remquanti} The following schemas are not \quci-valid:\\[1mm]
\begin{tabular}{llll}
(1) & $\ex x\no\varphi\imp\no\pt x\varphi$, & (2) & $\pt x\no\varphi\imp\no\ex x\varphi$, \\[1mm]
(3) & $\pt x\varphi\imp\no\ex x\no\varphi$, & (4) & $\ex x\varphi\imp\no\pt x\no\varphi$.\\[1mm]
\end{tabular}

\noindent
Indeed, let us consider the signature $\Sigma=\brac P$ where $P$ is an unary predicate symbol and consider the \quci-structure $\ff A$ where $A=\{a,b,c\}$  and $P^\ff A=\brac{ \ a \ , \ \emptyset \ , \ \{b,c\} \ }$. For any variable $x$ and any $e \in A$ let $[x\mapsto e] := \{s \in S(\ff A)  \ : \  s(x)=e\}$.
Then, in $\ff A$ we have
$$
\begin{array}{lll}
\|P(x)\|^\ff A
&=&
\brac
{
\ [x\mapsto a] \ ,
\ \emptyset \ ,
\ [x\mapsto b]\cup[x\mapsto c] \ 
}\\[2mm]
\mas{\pt xP(x)} 
&=&
\widehat{\ex x}([x\mapsto a]) \ = \
S(\ff A)
\\[2mm]
\menos{\no\pt x P(x)}
&=& \mas{\pt xP(x)} \ = \ S(\ff A)
\\[2mm]
\|\no P(x)\|^\ff A
&=&
\brac
{
\ \emptyset \ ,
\ [x\mapsto a] \ ,
\ [x\mapsto b]\cup[x\mapsto c] \ 
}\\[2mm]
\menos{\ex x\no P(x)}
&=&
\widehat{\pt x}([x\mapsto a]) \ = \ \emptyset.
\end{array}$$
Therefore, it is not the case that \, $\ff A\dodash\ex x\no\varphi\imp\no\pt x\varphi$, by Proposition~\ref{citad.imp}~(iii).
In a similar way, we can construct \quci-structures in order to show that (2), (3) and (4) are not valid in \quci. This shows that \quci\ differs from ${\bf LFI2^*}$, the first-order version of \ciore\ proposed in~\cite{lf79}. Indeed, the formulas $\no\pt x\varphi$ and $\ex x\no\varphi$ are equivalent in ${\bf LFI2^*}$, as well as the formulas $\no\ex x\varphi$ and $\pt x\no\varphi$.
\end{remark}

\

\noindent From the results stated in this section, we have:

\begin{theorem} [Soundness] \label{soundQmbC} 
Let $\Gamma \cup \{\varphi\}$ be a set of formulas over $\Gamma$.
If  $\Gamma\vdash_{\bf QCiore}\varphi$ \ then \ $\Gamma\dodash\varphi$.
\end{theorem}
\begin{dem} By induction on the complexity of  $\varphi$ and taking into account that all the axioms of \quci\ are valid (see Corollary~\ref{cim.sietee} and  Proposition~\ref{cim.diez}) and that all the  inference rules of \quci\ preserve validity (see  Propositions~\ref{cim.tres} and~\ref{cim.once}).
\end{dem}

\

In order to prove completeness, some well-known notions and results from classical first-order logic will be adapted to this framework. Recall that a {\em theory} is any set of sentences (that is, a set of formulas without free variables). A theory $T\subseteq\sent\Sigma$ is {\em \quci-consistent} if for every $\varphi\in\sent\Sigma$ at most two of  $\varphi$, $\no\varphi$, $\cons\varphi$  are derivable from $T$ in \quci\ (recalling from Section~\ref{s2} that $\sent\Sigma$ denotes the set of sentences over $\Sigma$).

\begin{proposition}\label{ciconsi1}
A theory $T$  is not \quci-consistent if and only if $T\vdash_{\bf QCiore}\varphi$ for every sentence $\varphi$.
\end{proposition}
\begin{dem} ($\Rightarrow$) Suppose that there exists a sentence $\psi$ such that \, $T\vdash_{\bf QCiore} \psi$, $T\vdash_{\bf QCiore} \no\psi$ and $T\vdash_{\bf QCiore} \cons\psi$; and let $\varphi$ be an arbitrary sentence. From {\bf (bc1)} and (MP) it follows that $T\vdash_{\bf QCiore} \varphi$.\\[1mm]
($\Leftarrow$) It is immediate from the definitions.
\end{dem}

\

\begin{corollary}\label{ciconsi2}
$T$ is \quci-consistent if and only if there exists a sentence that is not derivable from $T$ in $\quci$.
\end{corollary}
\begin{dem} It is a consequence of Proposition~\ref{ciconsi1}.
\end{dem}

\

\begin{proposition}\label{ciconsi3}
Let $T$ be a theory and let $\varphi$ be a sentence.
\begin{itemize}
\item[(i)]  $\varphi$ is not derivable from $T$ in \quci\ \, iff \, $T\cup\nono\varphi$ is \quci-consistent,
\item[(ii)] $\no\varphi$ is not derivable from $T$ in \quci\ \, iff \, $T\cup\sisi\varphi$ is \quci-consistent,
\item[(iii)] $\cons\varphi$ is not derivable from $T$ in \quci\ \, iff \, $T\cup\nosi\varphi$ is \quci-consistent.
\end{itemize}
\end{proposition}
\begin{dem}
It is consequence, respectively, of the following schemas of  \eledos-tautologies (which are therefore theorems of \quci): $((\no\varphi\wedge\cons\varphi)\imp\varphi)\imp\varphi$, $((\varphi\wedge\cons\varphi)\imp\no\varphi)\imp\no\varphi$ and $((\varphi\wedge\no\varphi)\imp\cons\varphi)\imp\cons\varphi$.
\end{dem}

\

We shall say that a theory $T$ is {\em maximal \quci-consistent } \ if \ $T$ is \quci-consistent and the only \quci-consistent theory which contains $T$ is $T$ itself.

\begin{proposition}\label{ciconsimax}
Let $T$ be a  maximal \quci-consistent theory. Then, for every sentences $\varphi$ and $\psi$:
\begin{itemize}
\item[(i)] $\varphi\in T$ iff $T\vdash_{\bf QCiore}\varphi$;
\item[(ii)] $\varphi\not\in T$ iff $\nono\varphi\subseteq T$;
\item[(iii)] $\no\varphi\not\in T$ iff $\sisi\varphi\subseteq T$;
\item[(iv)] $\cons\varphi\not\in T$ iff $\nosi\varphi\subseteq T$;
\item[(v)] $\varphi\wedge\psi\in T$ \ iff \ $\varphi\in T$ and $\psi\in T$;
\item[(vi)] $\varphi\vee\psi\in T$ \ iff \ $\varphi\in T$ or $\psi\in T$;
\item[(vii)] exactly two of $\varphi$, $\no\varphi$, $\cons\varphi$ are in $T$;
\item[(viii)] if $\pt x\varphi\in T$, then $\varphi(c/x)\in T$ for every $c\in \mm C$.
\end{itemize}
\end{proposition}
\begin{dem}
\noindent (i) 
\begin{align*}
T\vdash_{\bf QCiore}\varphi\\
\Longrightarrow \quad &
T\not\vdash_{\bf QCiore}\no\varphi
\ \mbox{ or } 
T\not\vdash_{\bf QCiore}\cons\varphi
&&\mbox{(by {\bf QCiore}-consistency)}
\\
\Longrightarrow \quad &
T\cup\sisi\varphi \mbox{ is \quci-consistent} \\
\phantom{\Longrightarrow \quad} &
\ \mbox{ or } \
T\cup\nosi\varphi \mbox{ is \quci-consistent} \
&&
\mbox{(by Proposition~\ref{ciconsi3}(ii)(iii))}
\\
\Longrightarrow \quad &
\varphi\in T
&&\mbox{(by maximality).}
\end{align*}

\noindent (ii), \ (iii), and \ (iv): They are consequence of Proposition~\ref{ciconsi3}.

\noindent (v) It follows from  {\bf (Ax4)}, {\bf (Ax5)} and {\bf (Ax3)}. 

\noindent (vi) ($\Rightarrow$) Assume that $\varphi\vee\psi\in T$.  If $\varphi\not\in T$ then 
$(\no\varphi\wedge\cons\varphi)\wedge(\varphi\vee\psi)\in T$, from items~(ii)  and~(v). Since $\vdash_{\bf QCiore}\big( 
(\no\varphi\wedge\cons\varphi)\wedge(\varphi\vee\psi) \big)\imp \psi$, it follows from item~(i) that $\psi\in T$.\\
($\Leftarrow$) It is a consequence of \acio6 and \acio7.

\noindent (vii) It follows from items (ii), (iii), (iv) and (v).

\noindent (viii) It is a consequence of \acio{12}.
\end{dem}

\

\noindent Given that $\quci$ is a Tarskian and finitary logic, the  Lindenbaum-\L os Lemma (see~\cite[Theorem~22.2]{woj:84} or~\cite[Theorem~2.2.6]{CarCon}) holds for this logic. In particular, this result holds when $\quci$ is restricted to sentences. Moreover, because of Corollary~\ref{ciconsi2}, the  Lindenbaum-\L os Lemma for the logic $\quci$ restricted to sentences is equivalent to the following:

\begin{theorem} [Lindenbaum-\L os Lemma] \label{cilinda}
Any \quci-consistent theory can be extended to a maximal \quci-consistent theory.
\end{theorem}

\

\noindent
A set of sentences $T$ is a {\em \quci-Henkin} theory \, iff \, for any formula $\varphi$ with at most one free variable $x$ there is a constant $c$ such that $\ex x\varphi\imp\varphi(c/x)\in T$. Then, by a simple adaptation of the well-known method of constants  introduced by L. Henkin (also using the weak version of the deduction meta-theorem, that is, Theorem \ref{WDMT} restricted to sentences)  we can prove the following:

\begin{theorem} \label{maxHen}
Every \quci-consistent theory $T$ over a signature $\Sigma$ has a maximal \quci-consistent and \quci-Henkin extension over a signature $\Sigma_C$ obtained from $\Sigma$ by adding a set $C$ of new individual constants. 
\end{theorem}

\

\begin{theorem}[Henkin Lemma]\label{hen}	
Let $T$ be a maximal \quci-consistent and {\bf QCiore}-Henkin theory. Then $T$ has a \quci-model.
\end{theorem}
\begin{dem} Let us define an \quci-structure $\ff A$ over $\Sigma$ as follows: 
\begin{itemize}
\item[-] the domain $A$  is the set of all closed terms over $\Sigma$,
\item[-] if  $c$ is a constant symbol of $\Sigma$ then  $c^\ff A:=c$,
\item[-] if $f$ is a $n$-ary function symbol of $\Sigma$, $f^\ff A:A^n\imp A$ is given by $f^\ff A(\vct a):=f(\vct a)$,
\item[-] if $R$ is a $n$-ary predicate symbol, $\triple R$ is defined as follows
\begin{align*}
(\vct a)\in\tas  R
&& \ \mbox{ iff } \ &&
\uno{R(\vct a)}\in T \\
(\vct a)\in\tenos  R
&& \ \mbox{ iff } \ &&
\cero{R(\vct a)}\in T \\
(\vct a)\in\tpor  R
&& \ \mbox{ iff } \ &&
\tres{R(\vct a)}\in T 
\end{align*}
\end{itemize}

\noindent By induction on the complexity of a  sentence $\varphi$ in  in $\sent\Sigma$ it will be  proven that
\begin{align*}
&(1)\quad
\tripc\varphi= 1
&&\ \mbox{ iff } \ &&
\uno\varphi\in T
\\
&(2)\quad
\tripc\varphi= 0
&&\ \mbox{ iff } \ &&
\cero\varphi\in T
\\
&(3)\quad
\tripc\varphi= \frac{1}{2}
&&\ \mbox{ iff } \ &&
\tres\varphi\in T.
\end{align*}

\paragraph{Base step:} Let $\varphi$ be the atomic formula $R(\vct a)$. Then,  $\tripc{R(\vct a)}= 1$  iff $(\vct a)\in\tas R$
iff (by definition) $\uno{R(\vct a)}\in T$. Hence, $\ff A$ fulfills condition (1). Similarly, we prove conditions (2) and (3).

\paragraph{(IH):} Suppose that the theorem holds for every formula $\psi$ whose complexity is less or equal than $n$. \\[2mm]
Observe that, in the inductive step (that is, when analyzing a formula $\varphi$ with complexity $n+1$), it is enough to prove the `only if' direction of (1), (2) and (3), because of Proposition~\ref{ciconsimax}.  Indeed, assuming that the `only if' direction of (1), (2) and (3) was proved by a given $\varphi$, the corresponding `if' direction  of (1), (2) and (3) can be proven as follows: if $\tripc{\varphi}= 1$ then $\tripc{\varphi} \neq 0$ and $\tripc{\varphi}\neq \frac{1}{2}$. Then, by the `only if' part of~(2) and~(3),  $\tres{\varphi}\notin T$ and $\cero{\varphi}\notin T$. By items~(v) and~(vii) of Proposition~\ref{ciconsimax}, $\uno{\varphi}\in T$. The cases $\tripc{\varphi}= 0$ and $\tripc{\varphi}=\frac{1}{2}$ are proved analogously.

Thus, in order to complete the proof let us prove, assuming (IH), the `only if' part of (1), (2) and (3) for a given formula $\varphi$ with complexity $n+1$. We have the following cases to consider:

\paragraph{Case a:} $\varphi$ is $\no\psi$

(a.1) If  $\uno{\no\psi}\in T$ then, by Theorem~\ref{esquema}(vi),  $\cero\psi\in T$ and so, by (HI), $\tripc{\psi}= 0$. From this $\tripc{\no\psi}= \no\tripc{\psi}=1$.

(a.2) If $\cero{\no\psi}\in T$ then, by {\bf (cf)} and Theorem~\ref{esquema}(vi), $\uno\psi\in T$ and so, by (HI),  $\tripc{\psi}= 1$. Then $\tripc{\no\psi}= \no\tripc{\psi}= 0$.

(a.3) If $\tres{\no\psi}\in T$ then, by {\bf (cf)}, $\tres\psi\in T$ and so, by (HI), 
$\tripc{\psi}=\frac{1}{2}$. From this $\tripc{\no\psi}= \no\tripc{\psi}=\frac{1}{2}$.

\paragraph{Case b:} $\varphi$ is $\cons\psi$

(b.1) If $\uno{\cons\psi}\in T$ then $\cons\psi\in T$  and so, by 
(HI)  and Proposition~\ref{ciconsimax}, $\tripc{\psi}\in\{1,0\}$. Then $\tripc{\cons\psi}= \cons\tripc{\psi}=1$.

(b.2) If $\cero{\cons\psi}\in T$ then  $\neg\cons\psi\in T$. By Theorem~ \ref{esquema}(iv),  $\tres\psi\in T$ and so, by (HI), $\tripc{\psi}=\frac{1}{2}$. This implies that $\tripc{\cons\psi}= \cons\tripc{\psi}=0$.

(b.3) Since $T$ is \quci-consistent and $\cons\cons\psi \in T$ (by Theorem~ \ref{esquema}(v)) it follows that $\cons\psi\wedge\no\cons\psi \notin T$. Then, this subcase holds vacuously.

\paragraph{Case c:} $\varphi$ is $\psi\wedge\xi$

(c.1) If $(\psi\wedge\xi) \wedge \cons(\psi\wedge\xi) \in T$ then, by Theorem~\ref{esquema}(ix), $\big((\psi\wedge(\xi\wedge\cons\xi)\big) \vee ((\psi\wedge\cons\psi)\wedge\xi)\big) \in T$. Hence, by (IH) and Proposition~\ref{ciconsimax}:  either
$\tripc\psi\in\{1,\frac{1}{2}\}$ and $\tripc\xi=1$, or  $\tripc\psi=1$ and $\tripc\xi\in\{1,\frac{1}{2}\}$. From this $\tripc{\psi\wedge\xi}=1$.

(c.2) If $\no(\psi\wedge\xi)\wedge\cons(\psi\wedge\xi) \in T$ then $(\no\psi\wedge\cons\psi)\vee(\no\xi\wedge\cons\xi)\in T$, by Theorem~\ref{esquema}(vii). Hence, by (IH) and Proposition~\ref{ciconsimax}: either
$\tripc\psi=0$ or  $\tripc\xi=0$. From this
$\tripc{\psi\wedge\xi}=\tripc{\psi}\wedge\tripc{\xi}=0$.

(c.3) If $(\psi\wedge\xi)\wedge\no(\psi\wedge\xi) \in T$ then $(\psi\wedge\no\psi)\wedge(\xi\wedge\no\xi)\in T$, by Theorem~\ref{esquema}(viii). By (IH) and Proposition~\ref{ciconsimax} it follows that $\tripc\psi=\tripc\xi=\frac{1}{2}$ and so $\tripc{\psi\wedge\xi}=\tripc{\psi}\wedge\tripc{\xi}=\frac{1}{2}$.

\paragraph{Case d:} $\varphi$ is $\psi\vee\xi$

The three subcases are treated as in {\bf Case c}, but now using Theorem~\ref{esquema}(x), Theorem~\ref{esquema}(xi) and Theorem~\ref{esquema}(xii), respectively.

\paragraph{Case e:} $\varphi$ is $\psi\to\xi$

The three subcases are treated as in {\bf Case c}, but now using Theorem~\ref{esquema}~(xiii), Theorem~\ref{esquema}(xiv) and Theorem~\ref{esquema}(xv), respectively.

\paragraph{Case f:} $\varphi$ is $\ex x\psi$

(f.1)  Suppose that $\uno{\ex x\psi}\in T$. Since $\ex x\psi\in T$ and $T$ is {\bf QCiore}-Henkin, there exists a constant $c$ such that  
$\psi(c/x)\in T$. By (IH), we have that $\tripc{\psi(c/x)}\in\{1, \frac{1}{2}\}$. 
Then, $\widehat{\pt x}(\menos\psi)=\emptyset$. Since $\cons\ex x\psi\in T$ then $\ex x\cons\psi\in T$, by ({\bf Ax13}). Using once again that  $T$ is a {\bf QCiore}-Henkin theory, we conclude that there exists a witness $d$ such that 
$\cons\psi(d/x)\in T$. By (IH), we have that $\tripc{\psi(d/x)}\in\{1, 0\}$. Then,  $\widehat{\pt x}(\por\psi)=\emptyset$ and so $\mas{\ex x\psi}=S(\ff A)$. From this  $\tripc{\ex x\psi}= 1$. 

\

(f.2)  If $\tres{\ex x\psi}\in T$ then, by Theorem~\ref{ciorequant}(iv), $ \pt x(\psi\wedge\no\psi)\in T$. By ({\bf Ax12}) we have that $\tres{\psi(c/x)}\in T$,  for all $c\in A$. Then, by (IH), $\por{\psi(c/x)}=S(\ff A)$  for all $c\in A$ and therefore $\por{\ex x\psi}=\widehat{\pt x}(\por{\psi})=S(\ff A)$. From this $\tripc{\ex x\psi}= \frac{1}{2}$. 

\

(f.3) If $\cero{\ex x\psi}\in T$ then, by Theorem~\ref{ciorequant}(ii),  $\pt x(\cero{\psi})\in T$. By  ({\bf Ax12})  we have that
$\cero{\psi(c/x)}\in T$, for all  $c\in A$. By (IH), $\menos{\psi(c/x)}=S(\ff A)$  for all $ c\in A$ and then $\menos{\ex x\psi}=\widehat{\pt x}(\menos{\psi})=S(\ff A)$.  From this $\tripc{\ex x\psi}= 0$.

\paragraph{Case g:} $\varphi$ is $\pt  x\psi$ 

(g.1) Suppose that $\pt x\psi\wedge\cons\pt x\psi\in T$. Since $\vdash_{\bf QCiore} (\pt x\psi\wedge\cons\pt x\psi)\imp(\ex x\psi\wedge\cons\ex x\psi)$ then $\ex x\psi\wedge\cons\ex x\psi\in T$. Since $T$ is {\bf QCiore}-Henkin, there exists a witness $c$ such that $\psi(c/x)\wedge\cons\psi(c/x)\in T$. By (IH), it follows that $\tripc{\psi(c/x)}= 1$ and so $\mas{\psi(c/x)}=\ese A$. Then $\widehat{\ex x}(\mas{\psi})=\ese A$.

Suppose now that there exists $s\in \widehat{\ex x}(\menos{\psi(x)})$. By (IH), there exists $d$ such that $\no\psi(d/x)\wedge\cons\psi(d/x)\in T$. Since $\pt x\psi\in T$  we conclude that $\psi(d/x)\wedge\no\psi(d/x)\wedge\cons\psi(d/x)\in T$,
which contradicts the fact that $T$ is \quci-consistent. Therefore $\widehat{\ex x}(\menos{\psi})=\emptyset$ and then $\mas{\pt x\psi}
= \widehat{\ex x}(\mas{\psi})-\widehat{\ex x}(\menos{\psi}) = \ese A$. From this $\tripc{\pt x\psi}= 1$.

\

(g.2) If $\cero{\pt x\psi}\in T$ then $\ex x(\cero{\psi})\in T$, by Theorem~\ref{ciorequant}(iii). Since $T$ is {\bf QCiore}-Henkin, there exists $c$ such that $\cero{\psi(c/x)}\in T$. By (IH), $\tripc{\psi(c/x)}= 0$ and so $\menos{\psi(c/x)}=\ese A$. Then, $\menos{\pt x\psi}=\widehat{\ex x}(\menos{\psi})=\ese A$.	From this $\tripc{\pt x\psi}= 0$.									

\

(g.3) If $\tres{\pt x\psi}\in T$ then $\pt x(\tres{\psi})\in T$, by Theorem~\ref{ciorequant}(v). Hence $\tres{\psi(c/x)}\in T$ for every $c\in A$. By (IH), we conclude that $\por{\psi(c/x)}=\ese A$ for all $c\in A$. Then,  $\por{\pt x\psi}=\widehat{\pt x}(\por{\psi})=\ese A$. 	From this $\tripc{\pt x\psi}= \frac{1}{2}$.
\end{dem}

\

\begin{theorem}[Completeness for sentences] \label{complQmbC} Let $\Gamma \cup \{\varphi\}$ be a set of sentences over $\Sigma$. Then  $\Gamma\dodash\varphi$ implies $\Gamma\vdash_{\bf QCiore} \varphi$.
\end{theorem}
\begin{dem} Suppose that  $\varphi$ is not derivable from $\Gamma$ in {\bf QCiore}. Then, by Proposition~\ref{ciconsi3}, we have that the theory $\Gamma \cup\{\no\varphi, \cons\varphi\}$ is  \quci-consistent. By  Theorem~\ref{maxHen}, there is a theory $T$ over a signature $\Sigma_C$ (obtained from $\Sigma$ by adding a set $C$ of new constants) which is maximal \quci-consistent  and \quci-Henkin, and such that $\Gamma \cup\{\no\varphi, \cons\varphi\} \subseteq T$. By
Theorem~\ref{hen}, there exists an \quci-model $\ff A$ over $\Sigma_C$ for $T$. Let $\ff {A'}$ be the reduct of $\ff A$ to $\Sigma$. Then, $\ff {A'}\Vdash_\quci \Gamma$ and $\ff {A'}\Vdash_\quci \no\varphi\wedge\cons\varphi$, and so $\ff {A'}\nVdash_\quci \varphi$. Therefore, 
we have that $\Gamma\nVdash_\quci\varphi$.
\end{dem}

\

\begin{remark} \label{obsquant}
Observe that completeness was proved only for sentences, while soundness was stated for formulas in general (recall Theorem~\ref{soundQmbC}). However, a completeness theorems for formulas in general (i.e., for formulas possibly having free variables) can be obtained from Theorem~\ref{complQmbC} by observing the following: for any formula $\psi$ let $(\forall)\psi$ be the universal closure of $\psi$, that is: if $\psi$ is a sentence then  $(\forall)\psi=\psi$, and if $\psi$ has exactly the variables $x_1,\ldots,x_n$ occurring free then  $(\forall)\psi=(\forall x_1)\cdots(\forall x_n)\psi$. If $\Gamma$ is a set of formulas then  $(\forall)\Gamma := \{(\forall)\psi  \ : \  \psi \in \Gamma\}$. Thus, it is easy to see that,   for every set $\Gamma \cup \{\varphi\}$ of formulas:  $\Gamma\vdash_{\bf QCiore} \varphi$ \ iff \   $(\forall)\Gamma\vdash_{\bf QCiore} (\forall)\varphi$, and $\Gamma \dodash \varphi$ \ iff \   $(\forall)\Gamma \dodash (\forall)\varphi$. From this, a general completeness result follows (see Corollary below).
\end{remark}

\begin{corollary} [Completeness] \label{complQmbC1} Let $\Gamma \cup \{\varphi\}$ be a set of formulas over $\Sigma$. Then  $\Gamma\dodash\varphi$ implies $\Gamma\vdash_{\bf QCiore} \varphi$.
\end{corollary}
\begin{dem}
It follows from Remark~\ref{obsquant} and by Theorem~\ref{complQmbC}.
\end{dem}

\begin{theorem} [Compactness] \label{com1}
A theory $T$ has a \quci-model iff every finite subset of $T$ has a \quci-model.
\end{theorem}
\begin{dem} ($\Rightarrow$): Immediate.\\[1mm] 
($\Leftarrow$): Since every sub-theory of $T$ has a \quci-model, by hypothesis, then it follows that every finite sub-theory of $T$ is \quci-consistent. Therefore  $T$ must be \quci-consistent, by finiteness of \quci, and so $T$ has a \quci-model,  by Theorem~\ref{hen}. 
\end{dem}

\

\begin{corollary} \label{compacidad} Let $\Gamma \cup \{\varphi\}$ be a set of sentences.
If $\Gamma\dodash\varphi$ then $\Gamma_0\dodash\varphi$, for some finite $\Gamma_0\subseteq \Gamma$.
\end{corollary}
\begin{dem} Suppose that for every finite $\Gamma_0\subseteq \Gamma$ it is not the case of $\Gamma_0\dodash\varphi$. By Proposition \ref{ciconsi3}(i) we have that $T_0=\Gamma_0\cup\{\no\varphi,\cons\varphi\}$ is  \quci-consistent and so, by Theorems~\ref{maxHen} and~\ref{hen},  $T_0$ has  a \quci-model, for every $\Gamma_0\subseteq \Gamma$ with $\Gamma_0$ finite. By Theorem \ref{com1}, we have that $\Gamma\cup\nono \varphi$ has a \quci-model. Therefore, it is not the case of $\Gamma\dodash\varphi$.
\end{dem}

\

\section{\quci-structures with equality}\label{igualdad}

Given a first-order signature $\Sigma$, we shall denote $\Sigma(\approx)$ the extension of  $\Sigma$ obtained by adding a binary predicate symbol  $\approx$ that we call {\em equality}. As suggested by its name,  $\approx$ will represent the identity predicate within the paraconsistent first-order logic \quci.

\begin{definition} \label{defeq} The logic {\bf QCiore}($\approx$) over the signature $\Sigma(\approx)$ is the extension of {\bf QCiore} over $\Sigma(\approx)$ obtained by adding the following schemas (as usual, we write $(t \approx t')$ instead of ${\approx}(t,t')$):

\begin{itemize}
\item[-] $\forall x(x\approx x)$
\item[-] $\forall x\forall y\big((x\approx y) \imp (\varphi\rightarrow \varphi[x\wr y])\big)$, if $y$ is a variable free for $x$ in $\varphi$
\end{itemize}
where $\varphi[x\wr y]$ denotes any formula obtained from $\varphi$ by replacing some, but not necessarily all, free occurrences of $x$ by $y$. 
\end{definition}

\

\noindent Let $\vec x=x_1,\ldots,x_n$ be a list of $n$ distinct  variables, and let $\vec t=t_1,\ldots,t_n$ be a list of terms. If $\varphi$ is a formula depending at most on the variables $x_1,\ldots x_n$ then $\varphi[{\vec t}/\vec x]$ will denote the formula obtained from $\varphi$ by simultaneously replacing the variable $x_i$ by $t_i$, for $1 \leq i \leq n$. 

\begin{proposition}\label{propterm} \ \ $\textbf{QCiore}(\approx)$ proves the following, for any terms $t$, $t'$ and $t''$:
\begin{enumerate}
\item $\vdash (t\approx t)$,
\item $\vdash (t\approx t')\rightarrow (t'\approx t)$,
\item $\vdash (t\approx t')\rightarrow((t'\approx t'')\rightarrow (t\approx t''))$,
\item $\vdash \big(\bigwedge_{i=1}^n (t_i \approx {t'}_i)\big) \to (f(t_1,\ldots,t_n) \approx f(t'_1,\ldots,t'_n))$ for every terms $t_1,t'_1,\ldots,t_n,t'_n$ and any function symbol $f$ of arity $n$,
\item $\vdash \big(\bigwedge_{i=1}^n (t_i \approx {t'}_i)\big) \to (\varphi[{\vec t}/\vec x] \to \varphi[{\vec t'}/\vec x])$ for every terms $t_1,t'_1,\ldots,t_n,t'_n$ and any formula without quantifiers $\varphi$ depending at most on the variables $x_1,\ldots x_n$.
\end{enumerate}
\end{proposition}
\begin{dem} The proof is identical to the one for classical first-order logic with equality. For a proof of  items 1-3 see, for instance, \cite[Proposition~2.23]{mendel}. Items 4-5 are an easy consequence of  axiom ({\bf Ax12}) and the second axiom for $\approx$ presented in Definition~\ref{defeq}. Indeed, let $\vec y=y_1,\ldots,y_n$ be a list of $n$ distinct  variables such that $x_i \neq y_j$ for any $i,j$. Since $\varphi$ has no quantifiers and it depends at most on the variables $x_1,\ldots x_n$, it follows that $\varphi(y_1/x_1)(y_2/x_2) \cdots(y_n/x_n)=\varphi[{\vec y}/\vec x]$. From this, and by applying $n$-times the second axiom in Definition~\ref{defeq} for $(x_1\approx y_1), (x_2\approx y_2), \ldots, (x_n\approx y_n)$, together with  axiom ({\bf Ax12}) and the transitivity of $\to$, it follows that $\vdash  \big(\bigwedge_{i=1}^n (x_i \approx y_i)\big) \to (\varphi \to \varphi[{\vec y}/\vec x])$. By Theorem~\ref{genrule}, $\vdash  \forall x_1\forall y_1\cdots\forall x_n\forall y_n\big(\big(\bigwedge_{i=1}^n (x_i \approx y_i)\big) \to (\varphi \to \varphi[{\vec y}/\vec x])\big)$. Finally the result follows by  axiom ({\bf Ax12}), using the fact that $\varphi$ has no quantifiers  and it depends at most on the variables $x_1,\ldots x_n$.
\end{dem}

\

\noindent For any set $A$ let $\Delta(A):=\{(a,b) \in A \times A  \ : \  a=b\}$ be the {\em diagonal} of $A$. It is worth noting that a \quci-structure $\ff A$ over $\Sigma(\approx)$ is a model of the sentence $\forall x(x\approx x)$, if and only if $\Delta(A) \subseteq {\approx}_{\oplus}^{\ff A}\,\cup\,{\approx}_{\odot}^{\ff A}$. This motivates the following:

\begin{definition} \label{LFI-normal} We say that a \quci-structure $\ff A$ over $\Sigma(\approx)$ has a partial equality (or a pragmatic equality, or a paraconsistent equality) if 
$${\approx}_{\oplus}^{\ff A}\,\cup\,{\approx}_{\odot}^{\ff A}\ = \Delta(A).$$
In this case we say that $\ff A$ is a \quci$(\approx)$-structure.
\end{definition}

\begin{proposition} \label{norm-estr}
Let  $\ff A$ be a \quci$(\approx)$-structure. Then $\ff A$  is a model of the logic {\bf QCiore}($\approx$). That is, $\ff A$ is a model of the axioms presented in Definition~\ref{defeq}.
\end{proposition}
\begin{dem}
Straightforward, from the definition of \quci$(\approx)$-structure.
\end{dem}

\

\begin{remark} Observe that a \quci$(\approx)$-structure  $\ff A$ validates a sentence $(t\approx t')$ iff the interpretations of $t$ and $t'$ in  $\ff A$ coincide. That is, the equality symbol has a {\em standard} or {\em normal} interpretation, according to the usual terminology adopted in first-order logic. However, a \quci$(\approx)$-structure  $\ff A$ is also allowed to validate a sentence of the form $\neg(t\approx t)$. This means that the equality symbol, even with a standard interpretation, is paraconsistent in the present approach. This is similar to the case for {\bf QmbC} and other quantified {\bf LFI}s with equality (see~\cite[Chapter~7]{CarCon}).
\end{remark}

\begin{definition} $\Gamma \cup \{\varphi\}$ be a set of formulas over $\Sigma(\approx)$. We say that {\em $\varphi$ is a $\quci(\approx)$-consequence of $\Gamma$}, denoted by 
 $\Gamma\Vdash_{\quci(\approx)}\varphi$ if, for every
\quci$(\approx)$-structure $\ff A$, we have that: $\ff A\Vdash_\quci\psi$ for every $\psi\in \Gamma$ implies that $\ff A\Vdash_\quci\varphi$.
\end{definition}

\begin{theorem}\label{lemequal} 
Let $T$ be a maximal ${\bf QCiore}(\approx)$-consistent and ${\bf QCiore}(\approx)$-Henkin theory. Then $T$ has a \quci$(\approx)$-model.
\end{theorem}
\begin{dem}  Let $T$ be a theory which is ${\bf QCiore}(\approx)$-consistent and  ${\bf QCiore}(\approx)$-Henkin. By Theorem~\ref{hen} we have that $T$ has a \quci-model $\ff A$ whose domain is the set $A$ of all the closed terms over the signature $\Sigma(\approx)$. Moreover,  
$t^{\ff A}(s)=t$ for every closed term $t$ and every  $s\in S(\ff A)$, and 
\begin{align*}
&(1)\quad
\tripc\varphi= 1
&&\ \mbox{ iff } \ &&
\uno\varphi\in T
\\
&(2)\quad
\tripc\varphi= 0
&&\ \mbox{ iff } \ &&
\cero\varphi\in T
\\
&(3)\quad
\tripc\varphi= \frac{1}{2}
&&\ \mbox{ iff } \ &&
\tres\varphi\in T
\end{align*}
for every  sentence $\varphi$ in $\sent{\Sigma(\approx)}$. Let $\equiv$ be the following relation on $A$: $t \equiv t'$ iff $(t \approx t') \in T$. By Proposition~\ref{propterm}(1)-(3) it follows that $\equiv$ is an equivalence relation.
Let $\bar a$ be the equivalence class of $a\in A$ determined by $\equiv$. Consider a new \quci-structure $\ff A^*$  defined as follows: 
\begin{itemize}
\item[-] the domain of $\ff A^*$ is the set $A^*=A/_\equiv$ of equivalence classes;
\item[-] if $R$ is a predicate symbol of arity $n$ then $\triple R$ is defined as follows:
\begin{align*}
(\bar a_1,\ldots,\bar a_n)\in  R^{\mathfrak A^*}_\oplus
&& \ \mbox{ iff } \ &&
\uno{R(a_1,\ldots,a_n)}\in T \\
(\bar a_1,\ldots,\bar a_n)\in R^{\mathfrak A^*}_\ominus
&& \ \mbox{ iff } \ &&
\cero{R(a_1,\ldots,a_n)}\in T \\
(\bar a_1,\ldots,\bar a_n)\in  R^{\mathfrak A^*}_\odot
&& \ \mbox{ iff } \ &&
\tres{R(a_1,\ldots,a_n)}\in T 
\end{align*}

\item[-] if $f$ is a function symbol of arity $n$ then
$f^{\ff A^*}(\bar a_1,\ldots,\bar a_n) = \overline{f(a_1,\ldots, a_n)}
$;
\item[-] if $c$ is a constant symbol then  \ $c^{\ff A^*} \ = \ \overline{c}$.
\end{itemize}

\noindent {\bf Fact:} The \quci-structure $\ff A^*$ is well-defined.\\
Indeed, suppose that $R$ is a predicate symbol of arity $n$ and $a_i\equiv a'_i$ for $1\leq i \leq n$ (recalling that, by definition, every $a_i$ and $a'_i$ is a closed term). Then, by Proposition~\ref{propterm}(5), $\uno{R(a_1,\ldots,a_n)}\in T$ iff $\uno{R(a'_1,\ldots,a'_n)}\in T$. Hence $\tas  R$ is well-defined. Analogously it can be proven that $\tenos  R$ and $\tpor  R$ are well-defined. Now, if $f$ is a function symbol of arity $n$ and $a_i\equiv a'_i$ for $1\leq i \leq n$ then, by Proposition~\ref{propterm}(4), $(f(a_1,\ldots, a_n)\approx f(a'_1,\ldots, a'_n)) \in T$, therefore $f^{\ff A^*}:(A^*)^n \to A^*$ is a well-defined function. This proves the {\bf Fact}. \\[1mm]
Observe that $\approx^{\ff A^*}$ is a paraconsistent equality. Indeed,  $(\bar a, \bar a') \in {\approx}_\oplus^{\ff A^*}\,\cup\,{\approx}_\odot^{\ff A^*}$ \ iff, by definition of the interpretation of the predicates in ${\ff A^*}$ and $\ff A$, \
$(a, a')\in {\approx}_\oplus^{\ff A}\,\cup\,{\approx}_\odot^{\ff A}$
\ iff \ $(a \approx a')\in T$ \ iff \ 
$\bar a=\bar a'$.
Thus, $\ff A^*$ is a $\bf LFI2(\approx)$-structure. 
Finally, $\ff A^*$ is a $\bf LFI2(\approx)$-model for $T$. Indeed, let $s\in S(\ff A)$ and $s^*\in S(\ff A^*)$ be defined by $s^*(x)=\overline{s(x)}$. Then, it can be proved by induction on the complexity of the formula $\varphi$ that 
\begin{center}
$\ff A\Vdash \varphi[s]$ 
\ iff \
$\ff A^*\Vdash \varphi[s^*]$.
\end{center}
Then, $\ff A^*\Vdash_\quci T$. That is, $\ff A^*$ is a \quci$(\approx)$-model of $T$.
\end{dem}

\

\begin{theorem} [Soundness and completeness of {\bf QCiore}$(\approx)$] \label{soundcompleq}
Let $\Gamma \cup \{\varphi\}$ be a set of sentences over $\Sigma(\approx)$. Then $\Gamma\vdash_{{\bf QCiore}(\approx)} \varphi$ if and only if  $\Gamma\Vdash_{\quci(\approx)}\varphi$.

\end{theorem}
\begin{dem}($\Rightarrow$): It is consequence of  Theorem~\ref{soundQmbC} and Proposition~\ref{norm-estr}.\\[1mm]
($\Leftarrow$): Suppose that  $\varphi$ is not derivable from $\Gamma$ in {\bf QCiore}$(\approx)$. Then, by Proposition~\ref{ciconsi3}, we have that the theory
$$T'=\Gamma \cup\{\no\varphi, \cons\varphi\} \cup Ax({\approx})$$
(where $Ax({\approx})$ is the set of axiom schemas for equality given in Definition~\ref{defeq})  is  \quci-consistent. By  Theorem~\ref{maxHen}, there is a theory $T$ over a signature $\Sigma({\approx})_C$ (obtained from $\Sigma({\approx})$ by adding a set $C$ of new constants) which is maximal \quci-consistent  and \quci-Henkin, and such that $T' \subseteq T$. Therefore, $T$ is  a maximal ${\bf QCiore}(\approx)$-consistent and ${\bf QCiore}(\approx)$-Henkin theory. Thus,  $T$ has a \quci$(\approx)$-model $\ff A$ over $\Sigma({\approx})_C$, by Theorem~\ref{lemequal}. Let $\ff {A'}$ be the reduct of $\ff A$ to $\Sigma({\approx})$. Then, $\ff {A'}$ is a \quci$(\approx)$-structure over $\Sigma({\approx})_C$ such that
$\ff {A'}\Vdash_{\quci} \Gamma$ and $\ff {A'}\Vdash_{\quci} \no\varphi\wedge\cons\varphi$, and so $\ff {A'}\nVdash_{\quci} \varphi$. Therefore, 
we have that $\Gamma\nVdash_{\quci(\approx)}\varphi$.
\end{dem}

\

\noindent As it was done for \quci, Theorem~\ref{soundcompleq} can be easily extended to formulas (possibly) having free variables.
It is not difficult to obtain a Compactness Theorem for \quci($\approx$).

\

\section{Some model-theoric results for \quci}

In what follows, we shall assume that in every signature we have at least the binary predicate symbol $\approx$ of equality. Morever, every \quci-structure will be assumed to be a \quci$({\approx})$ structure, in the sense of Definition~\ref{LFI-normal}. Thus, we will simply write $\Sigma$ instead of $\Sigma({\approx})$ to refer to a given first-order signature.
The notion of partial equality together with the Compactness Theorem allow us to make model-theoretic constructions analogous to the ones made for Tarskian structures. To begin with,  the following theorems can be easily obtained (the details are left to the reader).

\begin{theorem}
If a theory $T$ has arbitrarily large finite \quci$({\approx})$-models, then it has an infinite \quci$(\approx)$-model. 
\end{theorem}

\begin{theorem}
If a theory  $T$ over a signature $\Sigma$ has infinite \quci$({\approx})$-models, then it has infinite \quci$({\approx})$-models of any given cardinality $\geqslant\aleph_0 +|\rr L(\Sigma)|$, where $|\rr L(\Sigma)|$ denotes the cardinality of the first-order language $\rr L(\Sigma)$ generated by  $\Sigma$.
\end{theorem}

\noindent
Let $\ff A$ and $\ff B$ be two \quci-structures over $\Sigma$. We say that $\ff A$ and $\ff B$ are {\em elementarily equivalent}, denoted by $\ff A\equiv\ff B$, if for every sentence $\varphi$, if $\ff A\dodash \varphi$ then $\ff B\dodash\varphi$. 

\begin{proposition}
Assume that  $\ff A\equiv\ff B$. Then, for every sentence $\varphi$ it holds: $\ff A\dodash \varphi$ iff $\ff B\dodash\varphi$. 
\end{proposition}
\begin{dem} The `only if' part is immediate, by definition of $\equiv$. Now, suppose that $\varphi$ is a sentence such that $\ff A\not\dodash \varphi$. By Corollary~\ref{corosent} it follows that $\ff A\dodash \neg\varphi \wedge\cons\varphi$. By hypothesis, $\ff B\dodash \neg\varphi \wedge\cons\varphi$ and so, using again Corollary~\ref{corosent}, it follows that  $\ff B\not\dodash \varphi$. 
\end{dem}

\

\noindent $\ff A$ is said to be a {\em sub-\quci-structure} of $\ff B$, denoted by  $\ff A \ \subseteq \ \ff B$, if it is verified:  (1) $c^\ff A=c^\ff B$ for every constant symbol $c$, (2) $f^\ff A=f^\ff B|_{A^n}$ for every $n$-ary function symbol $f$; and (3) $R^\ff A_\#=R^\ff B_\#\cap A^n$ for every $\#\in\{\oplus,\odot,\ominus\}$ and for every $n$-ary predicate symbol $R$.

\begin{remark} \label{sbar} Suppose that $\ff A \ \subseteq \ \ff B$. If $s\in S(\ff A)$ then it induces a unique $\bar s\in S(\ff B)$ given by $\bar s(x)=s(x)$ for every $x \in \mm V$. The subset $\{ \bar s  \ : \  s \in S(\ff A)\}$ of $S(\ff B)$ will be denoted by $\bar S(\ff A)$. Note that $\bar{s}_x^a=\overline{s_x^a}$ for every $s\in S(\ff A)$, $x \in \mm V$ and $a \in A$. For every formula  $\varphi$ and $\#\in\{\oplus,\odot,\ominus\}$ let
$$\overline{\|\varphi\|^\ff A_\#} = \{\bar s   \ : \  s \in \|\varphi\|^\ff A_\#\}.$$
\end{remark}

\begin{definition} \label{elementary}
We say that $\ff A$ is a elementary sub-\quci-structure of $\ff B$, denoted by $\ff A\prec\ff B$, if   $\ff A\subseteq\ff B$ and for every formula  $\varphi$ it holds:
  
$$\overline{\|\varphi\|^\ff A_\#} = \|\varphi\|_\#^\ff B\cap \bar S(\ff A)$$

\noindent where $\#\in\{\oplus,\odot,\ominus\}$. 
Equivalently, $\ff A\prec\ff B$ iff for every formula $\varphi$ and every $s\in S(\ff A)$ we have that

\begin{center}
$\ff A\dodash \varphi[s]$ \ iff \ $\ff B\dodash\varphi[\bar s]$
\end{center}
\end{definition}

\noindent From the definition, it clearly follows that  $\ff A\prec\ff B$ implies $\ff A\equiv \ff B$.

\begin{lemma}[Tarski-like conditions] \label{Tcond}
Let $\ff A$ and $\ff B$ be two \quci-structures such that $\ff A\subseteq\ff B$. Then, the following conditions are equivalent:
\begin{itemize}
\item[(i)] $\ff A\prec\ff B$,
\item[(ii)] for all $s\in S(\ff A)$ the following conditions hold:
\begin{itemize}
\item[(TC1)] If $\bar s\in \mass{\ex x\varphi}B$, then there exist $a,b\in A$ such that $\overline{s_x^a}\not\in\menoss\varphi B$ and $\overline{s_x^b}\not\in\porr\varphi B$
\item[(TC2)] 
If $\bar s\not\in\porr{\ex x\varphi}B$, then there exists $a\in A$ such that $\overline{s_x^a}\not\in\porr\varphi B$
\item[(TC3)] If $\bar s\in\mass{\pt x\varphi}B$, then there exists $a\in A$ such that $\overline{s_x^a}\in\mass\psi B$
\item[(TC4)] If $\bar s\in\menoss{\pt x\varphi} B$, then there exists $a\in A$ such that $\overline{s_x^a}\in\menoss\varphi B$.

\end{itemize}
\end{itemize}
\end{lemma}
\begin{dem} $(ii) \Rightarrow (i)$: By induction on the complexity of the formula $\varphi$ it will be proven that 
\begin{center}
$\overline{\|\varphi\|_\#^\ff A}=\|\varphi\|_\#^\ff B\cap \bar S(\ff A)$.  
\end{center}
{\bf Base step:\ } $\varphi$ is an atomic formula $P(t_1,\ldots,t_n)$. Then, for each $\#\in\{\oplus,\ominus,\odot\}$: $s\in\|P(t_1,\ldots,t_n)\|_\#^\ff A$ iff $\big(t_1^\ff A[s],\ldots,t_n^\ff A[s]\big)\in P_\#^\ff A$ iff
$\big(t_1^\ff B[\bar s],\ldots,t_n^\ff B[\bar s]\big)\in P_\#^\ff B$ (since $\ff A\subseteq \ff B$) iff $\bar s\in\|P(t_1,\ldots,t_n)\|_\#^\ff B$.\\[1mm]
{\bf Inductive step:} Assume that $\overline{\|\psi\|_\#^\ff A}=\|\psi\|_\#^\ff B\cap \bar S(\ff A)$ for every $\psi$ with complexity $\leq k$. Let $\varphi$ with complexity $k+1$, and  let $s\in S(\ff A)$.\\[1mm]
(1) If $\varphi$ is one of $\no\psi$, $\circ\psi$, $\psi\wedge\chi$, $\psi\vee\xi$ or $\psi\imp\xi$, it is proved easily by using the induction hypothesis. \\[1mm]
(2) If $\varphi$ is $\ex x \psi$ then: \\[1mm]
(2.1) $s\in\mass{\ex x\psi}A$ only if $s \not\in \widehat{\pt x}(\menos\psi)\cup\widehat{\pt x}(\por\psi)$, by Proposition~\ref{charLFI2structure}, only if  there exist $a,b\in A$ such that $s_x^a\not\in\menos\psi $ and $s_x^b\not\in\porr\psi A$, by Definition~\ref{def-quant}, only if  there exist $a,b\in A$ such that $\overline{s_x^a}\not\in\menoss\psi B$ and $\overline{s_x^b}\not\in\porr\psi B$  (by I.H. and since $s_x^a,s_x^b\in S(\ff A)$), only if  there exist $a,b\in B$ such that $\bar{s}_x^a\not\in\menoss\psi B$ and $\bar{s}_x^b\not\in\porr\psi B$ (since $A\subseteq B$), only if $\bar s\in \mass{\ex x\psi}B$, by Proposition~\ref{charLFI2structure}.  Conversely: $\bar s\in \mass{\ex x\psi}B$ only if there exist $a,b\in A$ such that $\overline{s_x^a}\not\in\menoss\psi B$ and $\overline{s_x^b}\not\in\porr\psi B$ (by (TC1))  only if $s_x^a\not\in\menoss\psi A$ and $s_x^b\not\in\porr\psi A$ (by I.H. and since $s_x^a,s_x^b\in S(\ff A)$) only if $s\in\mass{\ex x\psi}A$.\\[1mm]
(2.2) Assume that $s\in \porr{\ex x\psi}A$, and suppose that $\bar s\not\in\porr{\ex x\psi}B$. By (TC2), there exists $a\in A$ such that $\overline{s_x^a}\not\in\porr\psi B$. By I.H., $s_x^a\not\in\porr\psi A$ for some $a\in A$, a contradiction. Then 
$\overline{\porr{\ex x\psi}A}\subseteq \porr{\ex x\psi}B\cap \bar S(\ff A)$.
Conversely: $\bar s\in\porr{\ex x\psi}B=\widehat{\pt x} (\porr\psi B)$ only if
$\bar{s}_x^a\in \porr\psi B$ for every $a\in B$ only if $\overline{s_x^a}\in \porr\psi B$ for every $a\in A$  (since $A\subseteq B$) only if
$s_x^a\in \porr\psi A$ for every $a\in A$ (by I.H.) only if
$s\in\widehat{\pt x} (\porr\psi A) = \porr{\ex x\psi}A$.\\[1mm]
(2.3) $s\in \menoss{\ex x\psi}A$ iff $s\not\in\mass{\ex x\psi}A\cup\porr{\ex x\psi}A$ iff
$\bar s\not\in\mass{\ex x\psi}B\cup\porr{\ex x\psi}B$ (by (2.1) and (2.2)) iff
$\bar s\in \menoss{\ex x\psi}B$. \\[1mm]
(3) If $\varphi$ is $\pt x \psi$ then: \\[1mm]
(3.1) $s\in\mass{\pt x\psi}A=\widehat{\ex x}(\mass\psi A)-\widehat{\ex x}(\menoss\psi A)$, only if:
$(1) \ s\in \widehat{\ex x}(\mass\psi A)$ \  and  \ $(2) \ s\not\in \widehat{\ex x}(\menoss\psi A)$, 
only if $s_x^a\not \in \menoss\psi A$, for every $a\in A$ (by (2)), only if $\overline{s_x^a}\not \in \menoss\psi B$, for every $a\in A$  (by I.H.),  only if $\bar s\not\in\menoss{\pt x\psi} B=\widehat{\ex x}(\menoss\psi B)$ (by (TC4)), only if $\bar s\in \widehat{\ex x}(\mass\psi B)$ (by (1) and  I.H.) only if
$\bar s\in\mass{\pt x\psi} B=\widehat{\ex x}(\mass\psi B)-\widehat{\ex x}(\menoss\psi B)$. Conversely: $\bar s\in\mass{\pt x\psi}B=\widehat{\ex x}(\mass\psi B)-\widehat{\ex x}(\menoss\psi B)$ only if:
$(1) \ \bar s\in \widehat{\ex x}(\mass\psi B)$ \ and \ $(2) \ \bar s\not \in \widehat{\ex x}(\menoss\psi B)$, only if there exists $a\in A$ such that $s_x^a\in\mass\psi A$  (by (TC3) and I.H.), only if $s\in \widehat{\ex x}(\mass\psi A)$, only if $s\not\in \widehat{\ex x}(\menoss\psi A)$  (by (2), $A\subseteq B$ and I.H.), only if $s\in\mass{\pt x\psi} A=\widehat{\ex x}(\mass\psi A)-\widehat{\ex x}(\menoss\psi A)$. \\[1mm]
(3.2) Suppose that $\bar s\not\in\porr{\pt x\psi}B$. Since $\porr{\pt x\psi}B = \porr{\ex x\psi}B$ then, by (TC2), there exists $a\in A$ such that $\overline{s_x^a}\not\in\porr\psi B$. By I.H., $s_x^a\not\in\porr\psi A$ for some $a\in A$, hence $s\not\in\widehat{\pt x}(\porr\psi A)=\porr{\pt x\psi}A$. Then  $\overline{\porr{\pt x\psi}A}\subseteq \porr{\pt x\psi}B\cap \bar S(\ff A)$. Conversely, $\bar s\in\porr{\pt x\psi}B$ only if $\bar s\in \widehat{\pt x}(\porr\psi B)$, only if $\bar{s}_x^a\in\porr\psi B$ for every $a\in B$, only if $\overline{s_x^a}\in\porr\psi B$ for every $a\in B$, only if $s_x^a\in\porr\psi A$ for every $a\in A$ (since $A\subseteq B$ and by I.H.), only if $s\in \widehat{\pt x}(\porr\psi A) = \porr{\pt x\psi}A$.\\[1mm]
(3.3) The proof for $\menoss{\pt x\psi}A$ is similar to the case~(2.3).
\end{dem}

\

\noindent
Let $\alpha$ be an ordinal number and let $\Sigma$ a first-order signature. A {\em chain of \quci-structures of length $\alpha$} is an increasing sequence of \quci-structures $\ff A_0\subseteq \ff A_1\subseteq \ldots\subseteq \ff A_\beta\subseteq \ldots$ (for $\beta < \alpha$)  over $\Sigma$.
The {\em union} of a chain as above is the \quci-structure $\ff A=\bigcup_{\beta < \alpha}\ff A_\beta$ over $\Sigma$ such that: (1)~the domain of $\ff A$ is the union $A=\bigcup_{\beta < \alpha} A_\beta$ of the respective domains of the  \quci-structures $\ff A_\beta$; (2)~$R_\#^{\ff A} \ := \ \bigcup_{\beta < \alpha}R_\#^{\ff A_\beta}$ for $\#\in\{\oplus,\odot,\ominus\}$ and for  every $n$-ary predicate symbol $R \in \mm P_n$; (3)~$f^{\ff A} \ := \ \bigcup_{\beta < \alpha}f^{\ff A_\beta}$ for every $n$-ary function symbol $f$; and (4)~$c^{\ff A} \ := \ c^{\ff A_0}$ for every constant symbol $c$. Note that $f^{\ff A}$ is a function, since the chain is increasing. Moreover, if $R\in\mm P_n$ then 
$R^{\ff A}=\langle R^{\ff A}_{ \oplus},R^{\ff A}_{ \ominus},R^{\ff A}_{ \odot}\rangle$
is a triple  over $A^n$. Then   $\ff A$ is, in fact, a \quci-structure over $\Sigma$.

\begin{remark} Recall that $\ff A$ is a \quci$(\approx)$-structure if \ ${\approx}_{\oplus}^{\ff A}\cup {\approx}_{\odot}^{\ff A}\ = \{(a,b) \in A \times A  \ : \  a=b\}$. Then, it is easy to see that the union of a chain of  \quci$(\approx)$-structures is a  \quci$(\approx)$-structure.
\end{remark}

\noindent
Using union of chains, it can be proved a version for \quci-structures of a well-known theorem obtained in 1957 by Tarski and Vaught, which generalizes the downward L\"owenheim-Skolem-Tarski theorem:

\begin{theorem} [Tarski and Vaught theorem for \quci-structures] Let $\Sigma$ be a  first-order signature.
Let $\ff B$ be a \quci-structure over $\Sigma$ of cardinality $\alpha$, and let $\beta$ be a cardinal such that $|\rr L(\Sigma)|\leqslant \beta\leqslant\alpha$. Then, there exists a \quci-structure $\ff A$ over $\Sigma$ of cardinality $\beta$ such that $\ff A\prec\ff B$.  
\end{theorem}
\begin{dem} Let $X\subseteq B$ of cardinality $\beta$.
We shall construct inductively a denumerable sequence of sets $A_0\subseteq A_1\subseteq \ldots\subseteq A_n\subset\ldots$ as follows:  $A_0=X\cup\{c^\ff B \ : \ c\in\mm C\}$. Now,  suppose that $A_n$ was already defined. Then:\\

(LS0) \ if $f\in\mm F_k$ then $f^\ff B(a_1,\ldots,a_k)\in A_{n+1}$ for all $(a_1,\ldots,a_k)\in A_n^k$.\\

\noindent Moreover, for every formula $\varphi$ and every $s\in S(\ff B)$ taking values in $A_n$:\\

(LS1) if $s\in\mass{\ex x\varphi}B$, we choose $a, b\in B$ such that $s_x^a\not\in\menoss\varphi B$ and $s_x^b\not\in\porr\varphi B$,\\[2mm]
\indent
(LS2) \ if $s\not\in\porr{\ex x\varphi}B$, we choose $c\in B$ such that $s_x^c\not \in \porr \varphi B$,\\[2mm]
\indent
(LS3) \ if $s\in\mass{\pt x\varphi}B$, we choose $d\in B$ such that $s_x^d\in\mass\varphi B$,\\[2mm]
\indent
(LS4) \ if $s\in\menoss{\pt x\varphi}B$, we choose $e\in B$ such that $s_x^e\in\menoss\varphi B$.\\

\noindent
Note that the existence of the elements $a,b,c,d,e$ is guaranteed by Proposition~\ref{charLFI2structure}. Besides, we add to the set $A_{n+1}$ the elements $\{a,b\}$ (or the elements $c$, $d$ or $e$) as appropriate. Let $\ff A$ be the \quci-structure with domain $A=\bigcup_{n \geq 0} A_n$ such that
\begin{itemize}
\item[-] $c^\ff A := c^\ff B$, for every $c \in \mm C$,
\item[-] if $f\in\mm F_k$, let $f^\ff A:=f^\ff B|_{A^k}$ (which is well-defined, by (LS0))
\item[-] if $R\in\mm P_k$, let $R_\#^\ff A:=R_\#^\ff B\cap A^k$ for every $\#\in\{\oplus,\odot,\ominus\}$.
\end{itemize} 

\noindent
{\bf Claim 1:} The domain $A$ has cardinality $\beta$ and $\ff A\subseteq\ff B$.\\
It is immediate from the definitions.\\[1mm]
{\bf Claim 2: } Given $s\in S(\ff A)$, there exists  
$s':\mm V\rightarrow A_k$ (for some $k\in\bb N$) such that 
$\|\varphi\|^\ff B(s) \ = \ \|\varphi\|^\ff B(s')$.\\
It follows by the fact that every formula has a finite number of free variables, and by Proposition~\ref{ci1.15}.\\[1mm]
{\bf Claim 3: }
$\ff A$ and $\ff B$ satisfy conditions  (TC1)-(TC4) of Lemma~\ref{Tcond} (Tarski-like conditions). \\
In order to prove this,  let $s\in S(\ff A)$.\\[1mm]
Condition (TC1): 
If $\bar s\in\mass{\ex x\varphi}B$ then,
by {\bf Claim 2} and (LS1), 
there are $a, b\in A_{k+1}\subseteq A$ 
such that 
$\overline{s_x^a}\not\in\menoss\varphi B$ and $\overline{s_x^b}\not\in\porr\varphi B$. \\[1mm]
Condition (TC2): 
If $\bar s\not\in\porr\varphi B$, by  {\bf Claim 2} and (LS2), 
there exists $a\in A_{k+1}\subseteq A$ such that $\overline{s_x^a}\not\in\porr\varphi B$.\\[1mm]
Condition (TC3): 
If $\bar s\in\mass{\pt x\varphi}B$,
by  {\bf Claim 2} and (LS3), 
there exists $a\in A_{k+1}\subseteq A$ 
such that $\overline{s_x^a}\in\mass\varphi B$.\\[1mm]
Condition (TC4): 
If $\bar s\in\menoss{\pt x\varphi}B$,
by  {\bf Claim 2} and (LS4), 
there exists $a\in A_{k+1}\subseteq A$ 
such that $\overline{s_x^a}\in\menoss\varphi B$.\\

\noindent As a consequence of Lemma~\ref{Tcond}, $\ff A\prec\ff B$.

\end{dem}

\

\noindent
An {\em elementary chain} is a chain of \quci-structures such that $\ff A_k\prec\ff A_n$ whenever $k<n$. The Elementary Chain Theorem can be adapted to  \quci-structures:

\begin{theorem}[Elementary Chain Theorem] Let $\brac{\ff A_n \ : \  n\in\bb N}$ be an elementary chain of \quci-structures. Then $\ff A_k\prec\bigcup_n\ff A_n$ for all $k\in\bb N$.
\end{theorem}
\begin{dem}  Let us prove that $\overline{\|\varphi\|_\#^{\ff A_k}}=\|\varphi\|_\#^{\ff A}\cap \bar S(\ff A_k)$ for every $k\in\bb N$  and $\#\in\{\oplus,\odot,\ominus\}$ by induction on the complexity of the formula $\varphi$. Thus, let $k\in\bb N$ and $s\in S(\ff A_k)$.\\
{\bf Basic step:}  
Let $\varphi$ be the atomic formula $P(t_1,\ldots,t_n)$, and $\#\in\{\oplus,\odot,\ominus\}$. 
Clearly $t^\ff A[\bar s]=t^{\ff A_k}[s]$ for each term $t$.
Hence
\begin{align*}
\|P(t_1,\ldots,t_n)\|^\ff A(\bar s)
&=P^\ff A\big(t_1^\ff A[\bar s],\ldots,t_n^\ff A[\bar s]\big)\\
&=P^{\ff A_k}\big(t_1^{\ff A_k}[s],\ldots,t_n^{\ff A_k}[s]\big)\\
&=\|P(t_1,\ldots,t_n)\|^{\ff A_k}(s)
\end{align*}
{\bf Inductive step:} Assume that $\overline{\|\varphi\|_\#^{\ff A_k}}=\|\varphi\|_\#^{\ff A}\cap \bar S(\ff A_k)$ for every $k\in\bb N$ , $\#\in\{\oplus,\odot,\ominus\}$ and $\psi$ with complexity $\leq n$. Let $\varphi$ with complexity $n+1$, let $k\in \bb N$ and $s\in S(\ff A_k)$.\\[1mm]
(1) If $\varphi$ is one of $\no\psi$, $\circ\psi$, $\psi\wedge\chi$, $\psi\vee\xi$ or $\psi\imp\xi$, it is proved easily by using the induction hypothesis. \\[1mm]
(2) If $\varphi$ is $\ex x \psi$, then: \\[1mm]
(2.1) If $\bar s\in\mass{\ex x\psi}A$ then  there exist $a,b\in A=\cup_n A_n$ such that $\overline{s_x^a}\not\in\menoss\psi A$ and $s_x^b\not\in\porr\psi A$. Thus, there exist $j\in\bb N$ and $a,b\in A_j$ such that $\overline{s_x^a}\not\in\menoss\psi A$ and $\overline{s_x^b}\not\in\porr\psi A$. If $j\leq k$ then there
exist $a,b\in A_k$ such that $s_x^a\not\in\menoss\psi {A_k}$ and $s_x^b\not\in\porr\psi {A_k}$, by  I.H., then $s\in\mass{\ex x\psi}{A_k}$. Otherwise, if  $k<j$ then $\ff A_k\prec\ff A_j$. Note that $s, s_x^a,s_x^b \in S(\ff A_j)$ up to names, in the sense of Remark~\ref{sbar}.
By I.H.,  $s_x^a\not\in\menoss\psi {A_j}$ and $s_x^b\not\in\porr\psi {A_j}$, then $s\in\mass{\ex x\psi}{A_j}$. Since
$\overline{\|\ex x\psi\|_\oplus^{\ff A_k}}=\|\ex x\psi\|_\oplus^{\ff A_j}\cap \bar S(\ff A_k)$ (where the notation of  Definition~\ref{elementary} is now used for $\ff A_k\prec\ff A_j$), it follows that $s\in\mass{\ex x\psi}{A_k}$. Conversely, if $s\in\mass{\ex x\psi}{A_k}$ then there  exist $a,b\in A_k$ such that $s_x^a\not\in\menoss\psi{A_k}$ and $s_x^b\not\in\porr\psi{A_k}$. Then, there exist $a,b\in A$ such that $\bar{s}_x^a\not\in\menoss\psi{A}$ and $\bar{s}_x^b\not\in\porr\psi{A}$, by I.H. and $A_k\subseteq A$. Hence $s\in\mass{\ex x\psi}A$.\\[1mm]
(2.2) Assume that $s\in \porr{\ex x\psi}{A_k}=\widehat{\pt x}(\porr\psi{A_k})$. Suppose that $\bar s\not\in\widehat{\pt x}(\porr\psi A)$. Then, there exists $a\in A$ such that $\bar{s}_x^a\not\in\porr\psi A$. Hence, there exists $j\in\bb N$ and $a\in A_j$ such that $\overline{s_x^a}\not\in\porr\psi A$. If $j\leq k$ then $s_x^a\not\in\porr\psi{A_k}$ for some $a\in A_k$, by I.H., a contradiction (since $s\in \widehat{\pt x}(\porr\psi{A_k})$). Otherwise, if $k<j$ then 
$\ff A_k\prec\ff A_j$, and so, by I.H., $s_x^a\not\in\porr\psi{A_j}$ for some $a\in A_j$. 
From this, $s\not\in \porr{\ex x\psi}{A_j}$ and so $s\not\in \porr{\ex x\psi}{A_k}$ (given that $\ff A_k\prec\ff A_j$), a contradiction. Therefore, $\bar{s}\in \widehat{\pt x}(\porr\psi A)=\porr{\ex x\psi}{A}$. Conversely, suppose that $\bar s\in \porr{\ex x\psi}{A}=\widehat{\pt x}(\porr\psi{A})$. Then $\bar{s}_x^a\in\porr\psi A$ for every $a\in A$ and so $\overline{s_x^a}\in\porr\psi A$ for every $a\in A_k$. By I.H.,  $s_x^a\in\porr\psi{A_k}$ for every $a\in A_k$, hence $s\in \widehat{\pt x}(\porr\psi{A_k}) = \porr{\ex x\psi}{A_k}$.\\[1mm]
(2.3) The proof for $\menoss{\ex x\psi}A$ follows from (2.1) and (2.2). \\[1mm]
(3) If $\varphi$ is $\pt x \psi$, then: \\[1mm]
(3.1) Suppose that $\bar s\in\mass{\pt x\varphi}A=\widehat{\ex x}(\mass \psi A)-\widehat{\ex x}(\menoss\psi A)$. Then, 
$$(\star) \hspace{1cm} s\not\in\widehat{\ex x}(\menoss\psi{A_j}) \  \mbox{ for every $j \geq k$.}$$ 
Indeed, let  $j\geq k$.  Since  $\bar s\not\in\widehat{\ex x}(\menoss\psi A)$ then $\overline{s_x^a}\not\in\menoss\psi A$, for every $a\in A$. In particular, $\overline{s_x^a}\not\in\menoss\psi A$ for every $a\in A_j$. By I.H. for $j$,  $s_x^a\not\in\menoss\psi {A_j}$, for every $a\in A_j$. Then $s\not\in\widehat{\ex x}(\menoss\psi{A_j})$, proving $(\star)$. Now, since $\bar s\in\widehat{\ex x}(\mass \psi A)$ there exists $a\in A$ such that $\overline{s_x^a}\in\mass\psi A$. Then, there exists $j\in\bb N$ and $a\in A_j$ such that $\overline{s_x^a}\in\mass\psi A$. If $j\leq k$ then $s_x^a\in\mass\psi{A_k}$, by I.H., and so  $s\in \widehat{\ex x}(\mass\psi{A_k})$. From this and $(\star)$ it follows that $s\in\mass{\pt x\psi}{A_k}$. Otherwise, if  $k<j$, $s_x^a\in\mass\psi{A_j}$ by I.H. for $j$. Then  $s\in\widehat{\ex x}(\mass\psi{A_j})$. From this and  $(\star)$, $s\in\mass{\pt x\psi}{A_j}$. Given that $\ff A_k\prec \ff A_j$, it follows that $s\in\mass{\pt x\psi}{A_k}$.\\[1mm] 
(3.2) The case for $\porr{\pt x\psi}{A_k}$ is treated analogously.\\[1mm]
(3.3) The proof for $\menoss{\pt x\psi}A$ follows easily from (3.1) and (3.2).
\end{dem}

\

\

\noindent
Let $\ff A$ be a \quci-structure for $\Sigma$. If $\emptyset\neq X\subseteq A$, we expand $\Sigma$ to a new signature  $\Sigma_X$ by adding
a new constant symbol $c_a$ for each $a\in X$. In this way, we can expand $\ff A$ to a  \quci-structure \ $\ff A_X=(\ff A, a)_{a\in X}$ for $\Sigma_X$ where $c_a^{\ff A_X}:=a$. 
Given the map $h:A\imp B$ where $B$ is the universe of the \quci-structure $\ff B$ for $\Sigma$, then $\ff B_{h(X)}=(\ff B,h(a))_{a\in X}$ is the expansion of $\ff B$ to a \quci-structure for $\Sigma_X$ where $c_a^{\ff B_{h(X)}}=h(a)$.

\

More generally, we can choose a sequence $\vec c=c_1,\ldots,c_n$ of new constant symbols that are not in $\Sigma$ and consider the new signature $\Sigma_{\vec c}$ which is the expansion of $\Sigma$ obtained by adding the constant symbols $\vec c$. If $\ff A$ is a \quci-structure over $\Sigma$, we denote by $(\ff A,\vec a)$ the expansion of $\ff A$ in the signature $\Sigma_{\vec c}$ in such a way that $c_i^{(\ff A,\vec a)}=a_i$. If $\vec x=x_1,\ldots,x_n$ is a sequence of $n$ distinct variables then, for every  $s\in S(\ff A)$ such that $s(x_i)=a_i$, we have that
\begin{itemize}
\item[-] $t^\ff A[s]=t(\vec c)^{(\ff A,\vec a)}$ for all term $t(\vec x)$ 
\item[-] $\|\varphi\|^\ff A(s) \ =  \ \|\varphi(\vec c)\|^{(\ff A,\vec a)}$
for every formula $\varphi(\vec x)$ in the language $\Sigma$. 
\end{itemize}

\

\noindent It is not difficult to prove the next useful lemma.

\begin{lemma}[Lemma on Constants] \label{oncon}  Let $T\subseteq Sent(\Sigma)$ be a theory and let $\varphi(\vec x)$ be a formula in $\Sigma$. Let $\vec c$ be a sequence of different constant symbols not appearing in $T$. Then,  $T\dodash\varphi(\vec c)$ \, iff \, $T\dodash\pt\vec x\varphi$.
\end{lemma}

\

\noindent
Let $\ff A$ and $\ff B$ be two \quci-structures. A function $\mapa hAB$ is an {\em elementary embedding} of $\ff A$ into $\ff B$,
noted $h:\ff A\prec\ff B$, iff for all formulas $\varphi$ and assignments $s\in S(\ff A)$, we have
\begin{center}
$s\in\|\varphi\|_\#^\ff A$ \ sii \ $h\circ s\in\|\varphi\|_\#^\ff B$
\end{center}
for each $\#\in\{\oplus,\odot,\ominus\}$.  If $h:\ff A\prec\ff B$ and $h$ is onto, we say that $h$ is an {\em isomorphism} between $\ff A$ and $\ff B$. In such a case, we shall say that $\ff A$ and $\ff B$ are {\em isomorphic} and we shall write $\ff A\simeq\ff B$. It is easy to prove that $\ff A\equiv\ff B$ whenever $\ff A\simeq\ff B$.

The {\em elementary diagram} of $\ff A$ is the theory 
${\rm Th}(\ff A_A)
\ := \ 
\{\varphi\in {\rm Sent}\big(\Sigma_A\big) \ : \ 
\ff A_A\dodash\varphi\}$. As in classical model theory, the notion of elementary diagram is a tool used to characterize  elementary embeddings.

\

\begin{proposition}\label{1.30b} Let $\ff A$ and $\ff B$ be two \quci-structures for $\Sigma$ and let $\mapa fAB$. The following conditions are equivalent:
\begin{itemize}
\item[(a)] $f$ is an elementary embedding of $\ff A$ into $\ff B$,
\item[(b)] there is an elementary extension $\ff C\succ\ff A$ and an isomorphism $g:\ff C\imp \ff B$ such that $f\subseteq g$.
\item[(c)] $(\ff B,f(a))_{a\in A}$ is a \quci-model of the elementary diagram of $\ff A$. 
\end{itemize}
\end{proposition}

\

\begin{remark}  Let $\Sigma$, $\Sigma_1$ and $\Sigma_2$ be signatures such that $\Sigma\subset\Sigma_1$ and $\Sigma\subset\Sigma_2$, if $\ff A$ and $\ff B$ are \quci-structures over $\Sigma_1$ and $\Sigma_2$, respectively, we shall write
$$\ff A\equiv_\Sigma\ff B$$
to indicate that the $\Sigma$-reducts of $\ff A$ and $\ff B$ are elementary equivalent and, we shall write 
$$
\embe f{\ff A}{\ff B}\Sigma
$$
to indicate that $f$ is an elementary embedding of $\ff A|\Sigma$ in $\ff B|\Sigma$.
\end{remark}

\

\begin{theorem}\label{ool5.1} Let $\ff A$ and $\ff B$ be two \quci-structures over $\Sigma$ and $\Sigma'$, respectively, such that $\Sigma\subset\Sigma'$. If $\ff A\equiv_\Sigma \ff B$, then there exists a \quci-structure $\ff C$ over $\Sigma'$ such that $\embe f{\ff A}{\ff C}\Sigma$ and $\embe g{\ff B}{\ff C}{\Sigma'}$.
\end{theorem}
\begin{dem} Let $\Sigma':=\Sigma_\ff A\cup\Sigma_\ff B$ (and suppose that $A\cap B=\emptyset$). Let us consider the theory

\begin{equation*}
T:=Th(\ff A_A)\cup Th(\ff B_B)\ \subseteq \ Sent(\Sigma').
\end{equation*}

\noindent Suppose that there exists $T_0\subseteq T$ (finite) which does not have models. By taking conjunctions, we have that there exist $\varphi(\vec a)\in Th(\ff A_A)\cap T_0$ and $\psi(\vec b)\in Th(\ff B_B)\cap T_0$ such that $\varphi\wedge\psi$ does not have \quci-models. Then, none \quci-structure over $\Sigma'_\ff B$ which satisfies $\psi(\vec b)$ can be extended to a \quci-structure with interpretations for the symbols $\vec a$ such that $\varphi(\vec a)$ is \quci-satisfied. Thus, 

\begin{align*}
&{Th(\ff B_B)}\dodash{\no\varphi(\vec a)\wedge\cons\varphi(\vec a)}\\
&\Rightarrow \ \
{Th(\ff B_B)}\dodash{\pt\vec x(\no\varphi(\vec x)\wedge\cons\varphi(\vec x))}
&&\mbox{(Lemma on Constants)}
\\ 
&\Rightarrow \ \
{\ff B_B}\dodash{\pt\vec x(\no\varphi(\vec x)\wedge\cons\varphi(\vec x))}
\\ 
&\Rightarrow \ \
{\ff A_A}\dodash{\pt\vec x(\no\varphi(\vec x)\wedge\cons\varphi(\vec x))}
\\ 
&\Rightarrow \ \
\ff A_A\dodash{\no\varphi(\vec a)\wedge\cons\varphi(\vec a)}.
\end{align*}

\noindent Therefore $\ff A_A\dodash{\varphi(\vec a)\wedge\no\varphi(\vec a)\wedge\cons\varphi(\vec a)}$. By Theorem \ref{compacidad}, we have that $T$ has a model $\ff C$. Let us define the functions $\mapa f{A}{C}$ and
$\mapa g{B}{C}$ such that $f(x)=x^\ff C$ (for $x \in A$) and $g(y)=y^\ff C$ (for $y \in B$). By Proposition \ref{1.30b},  $f$ and $g$ are elementary embeddings.
\end{dem}

\

\begin{theorem}[Robinson's Joint Consistency Theorem]\label{oo5.2} 
Let $\Sigma_1$ and $\Sigma_2$ be two signatures such that $\Sigma=\Sigma_1\cap\Sigma_2$. Suppose that
\begin{itemize}
\item[(i)] $T\subseteq Sent(\Sigma)$ is a maximal \quci-consistent theory,
\item[(ii)] $T_1\subseteq Sent(\Sigma_1)$ and $T_2\subseteq Sent(\Sigma_2)$ extend to $T$, and
\item[(iii)] both $T_1$ and $T_2$ have a \quci-model. 
\end{itemize}
Then, $T_1\cup T_2$ has a \quci-model.
\end{theorem}
\begin{dem} Let $\ff A_0$ be a \quci-model of $T_1$ and $\ff B_0$ be a \quci-model of $T_2$. Then, $\ff A_0$ and $\ff B_0$ are models of $T$ and, since $T$ is maximal \quci-consistent, we have that \ $\ff A_0\equiv_\Sigma\ff B_0$.
By Theorem~\ref{ool5.1} and Proposition~\ref{1.30b}, there exists a \quci-structure $\ff B_1$ such that $\embe{f_0}{\ff A_0}{\ff B_1}{\Sigma}$. Since $(\ff A_0,A_0)\equiv_{\Sigma_{A_0}}(\ff B_1,f_0(A_0))$, again by Theorem \ref{ool5.1} and Proposition \ref{1.30b},  there exists a \quci-structure $\ff A_1$ $\embe{g_0}{\ff B_1}{\ff A_1}{\Sigma}$. Repeating this process, we can build a tower like the following:

\begin{center}
\begin{tikzpicture}[description/.style={fill=white,inner sep=2pt}]
\matrix (m) [matrix of math nodes, row sep=3em,
column sep=2.5em, text height=1.5ex, text depth=0.25ex]
{ \ff A_0 &  \ff A_1& \ff A_2&\cdots\\
\ff B_0 &\ff B_1 &\ff B_2&\cdots\\ };
\path[->,font=\scriptsize]
(m-1-1) edge node[auto] {$f_0$} (m-2-2)
edge node[auto] {$\preceq$} (m-1-2)
(m-2-1) edge node[auto] {$\preceq$}  (m-2-2)
(m-2-2) edge node[right] {$g_0$} (m-1-2)
(m-1-2) edge node[auto] {$f_1$} (m-2-3)
edge node[auto] {$\preceq$} (m-1-3)
(m-2-2) edge node[auto] {$\preceq$}  (m-2-3)
(m-2-3) edge node[right] {$g_1$} (m-1-3)
(m-1-3) edge (m-1-4)
(m-2-3) edge (m-2-4);
\end{tikzpicture}
\end{center}
where $\embe{f_k}{\ff A_k}{\ff B_{k+1}}\Sigma$ \ and \ $\embe{g_k}{\ff B_{k+1}}{\ff A_{k+1}}\Sigma$ are such that 

$$f_k\subseteq g_k^{-1}\subseteq f_{k+1}.$$

\noindent This is consequence of the characterization for embeddings given  by means of elementary diagrams.
By the Elementary Chain Theorem, $\ff A=\bigcup\ff A_k$ is a \quci-model of $T_1$ and $\ff B=\bigcup\ff B_k$ is a \quci-model of $T_2$. Besides, $\bigcup f_k$ is an isomorphism of $\ff A|\Sigma$ onto $\ff B|\Sigma$. 
Then, $\ff B$ is isomorphic to a \quci-structure $\ff B'$ such that $\ff A|\Sigma=\ff B'|\Sigma$. Let $\ff C$ be the \quci-structure on $\Sigma_1\cup\Sigma_2$ such that
\begin{itemize}
\item[-] $|\ff C|=|\ff A|$.

\item[-] 
If $c$ is a constant symbol in $\Sigma_2$, then $c^{\ff C}:=c^{\ff B'}$.

If $c$ is a constant symbol in $\Sigma_1-\Sigma_2$, then $c^{\ff C}:=c^{\ff A}$.

\item[-]
If $f$ is a function symbol in $\Sigma_2$, $f^\bb C=f^{\ff B'}$.

If $f$ is a function symbol in $\Sigma_1-\Sigma_2$, $f^\bb C=f^{\ff A}$.

\item[-] If $R$ is a predicate symbol in $\Sigma_2$, $R^\ff C=\triplee R{B'}$.
 
If $R$ is a predicate symbol in $\Sigma_1-\Sigma_2$, $R^\ff C:=\triplee RA$.
\end{itemize}
Then, $\ff C$ is a model of $T_1\cup T_2$.
\end{dem}

\

\begin{theorem}[Amalgamation] \label{oot5.3}
Every diagram
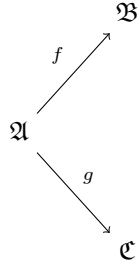
\begin{center}
\begin{tikzpicture}[description/.style={fill=white,inner sep=2pt}]
\matrix (m) [matrix of math nodes, row sep=3em,
column sep=2.5em, text height=1.5ex, text depth=0.25ex]
{ & \ff B\\
\ff A&\\
 &\ff C\\ };
\path[->,font=\scriptsize]
(m-2-1) edge node[auto] {$f$} (m-1-2)
edge node[auto] {$g$} (m-3-2);
\end{tikzpicture}
\end{center}
of elementary embeddings between \quci-structures on $\Sigma$, can be completed to a commutative diagram
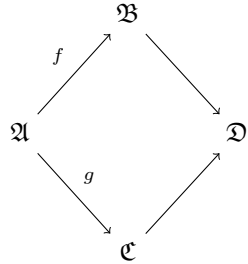
\begin{center}
\begin{tikzpicture}[description/.style={fill=white,inner sep=2pt}]
\matrix (m) [matrix of math nodes, row sep=3em,
column sep=2.5em, text height=1.5ex, text depth=0.25ex]
{ & \ff B&\\
\ff A&&\ff D\\
 &\ff C&\\ };
\path[->,font=\scriptsize]
(m-2-1) edge node[auto] {$f$} (m-1-2)
edge node[auto] {$g$} (m-3-2)
(m-1-2) edge (m-2-3)
(m-3-2) edge (m-2-3)
;
\end{tikzpicture}
\end{center}
of elementary embeddings between \quci-structures on $\Sigma$.
\end{theorem}
\begin{dem} Let us restate the elementary diagram $Th(\ff B, B)$ in the signature $\Sigma_1=\Sigma(A\cup(B-f(A)))$. In this signature, 
$c_a^\ff B=f(a)$, if $a\in A$. Analogously, we state the elementary diagram $Th(\ff C, C)$ in the signature $\Sigma_1=\Sigma(A\cup(C-g(A)))$ and, in this signature, $c_a^\ff C=g(a)$ if $a\in A$. Then, $\Sigma_1\cap\Sigma_2=\Sigma(A)$. \\
In $\Sigma_A$ we have a maximal \quci-consistent theory $Th(\ff A,A)$, and $Th(\ff B,B)$ and $Th(\ff C,C)$ are extensions of it in $\Sigma_1$ and $\Sigma_2$, respectively. By Theorem \ref{oo5.2}, $Th(\ff B,B)\cup Th(\ff C,C)$ has a model $\ff D$. \\
Now, define the functions 
\begin{center}
$
\begin{array}{rrcl}
f^*:&\ff B&\rightarrow&\ff D\\
&x&\mapsto&c_x^\ff D
\end{array}
$
\quad and \quad 
$
\begin{array}{rrcl}
g^*:&\ff C&\rightarrow&\ff D\\
&x&\mapsto&c_x^\ff D
\end{array}
$
\end{center}
Since $(\ff D,f^*(B))\dodash Th(\ff B,B)$ and $(\ff D,g^*(C))\dodash Th(\ff C,C)$ and by Proposition \ref{1.30b} we have that $f^*$ and $g^*$ are elementary embeddings.\\
Let $a\in\ff A$, then $f^*(f(a))=f^*(c_a^\ff B)=c_a^\ff D=g^*(c_a^\ff C)=g^*(g(a))$. Therefore,  $f^* \circ f=g^* \circ g$ (i.e. the diagram commutes).
\end{dem}

\

\begin{theorem}[Craig's Interpolation Theorem]\label{ch2.2.20} Let $\varphi$, $\psi$ be two sentences such that $\varphi\lfdash\psi$. Then, there exists a sentence $\theta$ such that 
\begin{itemize}
\item[(i)] $\varphi\dodash\theta$ and $\theta\dodash\psi$.
\item[(ii)] Every symbol of relation (excluding the identity), function or constant that occurs in $\theta$, it also occurs in both $\varphi$ and $\psi$.
\end{itemize}
\end{theorem}

\

\noindent In order to prove the Craig's Interpolation Theorem, we shall prove first the following technical results. \\
Let $T\subseteq Sent(\Sigma_1)$ and $U\subseteq Sent(\Sigma_2)$. We say that a sentence $\theta\in\Sigma_0$ {\em separates} the theories $T$ and $U$ if \ $T\dodash\theta \ \mbox{and} \ U\dodash {\no\theta\wedge\cons\theta}$. We say that $T$ are $U$ {\em inseparable} if no sentence $\theta\in Sent(\Sigma_0')$ separates them.

\

\begin{lemma}\label{lemma} If $T$ and $U$ are inseparable, then both are \quci-consistent.
\end{lemma}
\begin{dem} It is routine.
\end{dem}

\begin{lemma}\label{lemaA} If $T\cup\{\varphi\}$ and $U$ are inseparable then: either $T\cup\sisi\varphi$ and $U$ are inseparable, or $T\cup\nosi\varphi$ and $U$ are inseparable.
\end{lemma}
\begin{dem} Suppose that $T\cup\sisi\varphi$ and $U$ are separable, and that $T\cup\nosi\varphi$ and $U$ are separable. Then, there are $\theta_1$ and $\theta_2$ such that
 \begin{equation*}
 T\cup\sisi\varphi\dodash\theta_1, \quad U\dodash\cero{\theta_1}
 \end{equation*}
  \begin{equation*}
 T\cup\nosi\varphi\dodash\theta_2, \quad U\dodash\cero{\theta_2}
 \end{equation*}
Since 
$\dodash\varphi\imp(\varphi\wedge\cons\varphi)\vee(\varphi\wedge\no\varphi)$
\ and \
 $\dodash(\sneg\alpha\wedge\sneg\beta)\imp\sneg(\alpha\vee\beta)$ it follows that
\begin{equation*}
 T\cup\{\varphi\}\dodash\theta_1\vee\theta_2, \quad U\dodash\cero{(\theta_1\vee\theta_2)}.
 \end{equation*}
 
 \noindent
From this,  $T\cup\{\varphi\}$ and $U$ are separable.
\end{dem}

\

\begin{lemma}\label{lemaB} If $T\cup\{\no\varphi\}$ and $U$ are inseparable then: either  $T\cup\nono\varphi$ and $U$ are inseparable, or $T\cup\nosi\varphi$ and $U$ are inseparable.
\end{lemma} 
\begin{dem} Suppose that $T\cup\nono\varphi$ and $U$ are separable, and that  $T\cup\nosi\varphi$ and $U$ are separable.
Then, there are  $\theta_1$ and $\theta_2$ such that
 \begin{equation*}
 T\cup\nono\varphi\dodash\theta_1, \quad U\dodash\cero{\theta_1}
 \end{equation*}
  \begin{equation*}
 T\cup\nosi\varphi\dodash\theta_2, \quad U\dodash\cero{\theta_2}
 \end{equation*}
Since $\dodash\no\varphi\imp((\no\varphi\wedge\cons\varphi)\vee(\varphi\wedge\no\varphi))$ it follows that
  \begin{equation*}
 T\cup\{\no\varphi\}\dodash\theta_1\vee\theta_2, \quad U\dodash\cero{(\theta_1\vee\theta_2)}.
 \end{equation*}

 \noindent
From this,  $T\cup\{\no\varphi\}$ and $U$ are separable.
\end{dem}

\

\noindent
Analogously, we prove the next

\begin{lemma}\label{lemaC} If $T\cup\{\cons\varphi\}$ and $U$ are inseparable then: either $T\cup\sisi\varphi$ and  $U$ are inseparable or $T\cup\nono\varphi$ and $U$ are inseparable.
\end{lemma} 

\begin{dem} (of the Craig's Interpolation Theorem)
Suppose that there does not exist a Craig interpolant $\theta$ of $\varphi$ and $\psi$. We shall prove that it cannot happen that $\varphi\dodash\psi$. For this, we shall construct a \quci-model of $\varphi\wedge(\no\psi\wedge\cons\psi)$. 
Let $\Sigma_1$ the signature with all the symbols that occur in $\varphi$ and $\Sigma_2$ the signature with all the symbols that occur in  $\psi$. Now, consider the signature

$$\Sigma_0:=\Sigma_1\cap\Sigma_2, \quad \Sigma:=\Sigma_1\cup\Sigma_2.$$

Let $\Sigma'$ be an extension of $\Sigma$ that is obtained by adding a countable set $C$ of new constant symbols and let

\begin{equation*}
\Sigma'_0=\Sigma_0\cup C, \quad 
\Sigma'_1=\Sigma_1\cup C, \quad 
\Sigma'_2=\Sigma_2\cup C. 
\end{equation*}

\paragraph{(1)  $\{\varphi\}$ and $\{\no\psi,\cons\psi\}$ are inseparable.} \

\noindent
Indeed, suppose that there exists $\theta(\vec c)\in Sent(\Sigma_0')$, where $\vec c$ is a tuple of new constants  $C$ and such that 
$\varphi\dodash\theta(\vec c)$ and $\{\no\psi,\cons\psi\}\dodash\no\theta(\vec c)\wedge\cons\theta(\vec c)$. Then,  $\theta(\vec c)\dodash\psi$ and, by the Lemma on Constants, we have that $\varphi\dodash\pt \vec x\theta(\vec x)$ and $\pt \vec x\theta(\vec x)\dodash\psi$. Therefore, $\varphi$ and $\psi$ have a Craig interpolant.

\
Let $\varphi_0,\varphi_1,\varphi_2,\ldots$ be an enumeration of all the sentences in $\Sigma_1$ and let $\psi_0,\psi_1,\psi_2,\ldots$ be an enumeration of all the sentences in $\Sigma_2$. Then, we shall construct two increasing sequences of theories in the following way

\begin{align*}
\{\varphi\}&= \  T_0\subseteq T_1\subset T_2\subseteq\ldots,\\
\nono\psi&= \ U_0\subseteq U_1\subseteq U_2\subseteq\ldots,
\end{align*}
in $Sent(\Sigma'_1)$ and $Sent(\Sigma'_2)$ respectively, such that:
\begin{itemize}
\item[-] if $T_m\cup\{\varphi_m\}$ and $U_m$ are inseparable, then $\varphi_m\in T_{m+1}$,

\item[-] if $U_m\cup\{\psi_m\}$ and $T_{m+1}$ are inseparable, then $\psi_m\in U_{m+1}$.
\end{itemize}

\noindent Now, let $T_\omega=\bigcup T_m$ \ and \  $U_\omega=\bigcup U_m$.

\paragraph{(2) \ $T_\omega$ and $U_\omega$ are inseparable.} \

\noindent
Indeed, if there is $\theta\in\sent{\Sigma_0'}$ such that $T_\omega\dodash\theta$ and $U_\omega\dodash\cero\theta$. By Theorem~\ref{compacidad}, there exist  finite sets $\Delta\subseteq T_\omega$ and $\Gamma\subseteq U_\omega$ such that $\Delta\dodash\theta$ and $\Gamma\dodash\cero\theta$. Then, there is $k\in\omega$ such that
\begin{equation*}
T_k\dodash\theta, \quad U_k\dodash\cero\theta.
\end{equation*}

\noindent
Then $T_\omega$ and $U_\omega$ are both \quci-consistent, by Lema \ref{lemma}.

\paragraph{(3)  $T_\omega$ is a maximal \quci-consistent theory in $\Sigma_1'$, and $U_\omega$ is a maximal \quci-consistent theory in $\Sigma_2'$.} 

\

\noindent
Suppose that $\varphi_m\not\in T_\omega$, $\no\varphi_m\not\in T_\omega$ and $\cons\varphi_m\not\in T_\omega$. Then, $T\cup\nosi{\varphi_m}$, $T\cup\sisi{\varphi_m}$ and $T\cup\nono{\varphi_m}$, are separable from $U_m$. Then, there are $\theta_1,\theta_2,\theta_3\in Sent(\Sigma'_0)$ such that

\begin{center}
$
{T_\omega}
\dodash
{\varphi_m\wedge\no\varphi_m\imp\theta_1}$, \quad  
${U_\omega}
\dodash
{\no\theta_1\wedge\cons\theta_1}$,

\

${T_\omega}\dodash{\varphi_m\wedge\cons\varphi_m\imp\theta_2}$, \quad  ${U_\omega}\dodash{\no\theta_2\wedge\cons\theta_2}$,

\

${T_\omega}\dodash{\no\varphi_m\wedge\cons\varphi_m\imp\theta_3}$, \quad  ${U_\omega}\dodash{\no\theta_3\wedge\cons\theta_3}$.
\end{center}
Then, since $\dodash(\varphi_m\wedge\no\varphi_m)\vee(\varphi_m\wedge\cons\varphi_m)\vee(\no\varphi_m\wedge\cons\varphi_m)$ 
and
$\dodash(\sneg\theta_1\wedge\sneg\theta_2\wedge\sneg\theta_3)\imp\sneg(\theta_1\vee\theta_2\vee\theta_3)$
it follows that
\begin{center}
${T_\omega}\dodash{\theta_1\vee\theta_2\vee\theta_3}$ \ and \
${U_\omega}\dodash{\no(\theta_1\vee\theta_2\vee\theta_3)\wedge\cons(\theta_1\vee\theta_2\vee\theta_3)}$,
\end{center}
and this contradicts the inseparability of $T_\omega$ and $U_\omega$. Then, $T_\omega$ is a maximal \quci-consistent theory in $\Sigma'_1$. 
If $\varphi_m\in T_\omega$, then $\varphi_m\in T_m$. Now, $\varphi_r=\no\varphi_m$ and $\varphi_s=\cons\varphi_s$.
By Lemma \ref{lemaA} and the \ciore-consistency of $T_\omega$, we can assert that $\varphi_r\in T_\omega$ or $\varphi_s\in T_\omega$.

If $\varphi_r=\no\varphi_m\in T_\omega$, then $\varphi_r\in T_r$. If $\varphi_s=\cons\varphi_s$, by Lemma \ref{lemaB} and the \ciore-consistency of $T_\omega$, we have that $\varphi_m\in T_\omega$ or $\varphi_s\in T_\omega$.

If $\varphi_s=\cons\varphi_m\in T_\omega$, then $\varphi_s\in T_s$. If $\varphi_r=\no\varphi_s$, by Lemma \ref{lemaC} and the \ciore-consistency of $T_\omega$, we have that $\varphi_m\in T_\omega$ or $\varphi_r\in T_\omega$.

In a similar way, we prove the maximality of $U_\omega$.

\paragraph{(4) \ $T_\omega\cap U_\omega$ is a maximal \quci-consistent theory in $\Sigma'_0$. \ } 

\ 

\noindent
Let $\sigma$ be a sentence of $\Sigma'_0$. By (3) and the inseparability, we have that either 
\begin{center}
$\nosi\sigma\subseteq T_\omega\cap U_\omega$ \ { or } \
$\sisi\sigma\subseteq T_\omega\cap U_\omega$ \ { or } \
$\nono\sigma\subseteq T_\omega\cap U_\omega$.
\end{center}

\noindent
Then, by the Robinson's Joint Consistency Theorem, $T_\omega\cup U_\omega$ has a model. Therefore, $\{\varphi,\no\psi,\cons\psi\}$ has a model and then it cannot be the case that $\varphi\dodash\psi$.
\end{dem}

\

\section{Conclusions}

In the present paper we have developed a model theory framework that was used successfully for studying a first-order version of the paraconsistent 3-valued logic \ciore\ called \quci. This model theory strongly relies on the notion of partial structures presented in Section~\ref{s3}. Based on it, it was shown in which manner classical  model theory can be adapted to this (more general) setting. In particular, important classical theorems of classical model theory were obtained in this framework,  such as: {\em Tarski and Vaught Theorem} (generalizing the {\em Downward L\"owenheim-Skolem-Tarski Theorem}),
{\em Robinson's Joint Consistency Theorem}, {\em Amalgamation Theorem} and {\em Craig's Interpolation Theorem}.
It is worth mentioning that the concept of pragmatic structure can be used to study other 3-valued {\bf LFI}s. For instance,  all the results obtained by  D'Ottaviano in \cite{it127, dot:85a, dot:85b,dot:87} for the first-order version of \dacdot\ can be recast in this framework.

An interesting task to be done is trying to adapt this framework for studying a first-order version of each of the 8K 3-valued {\bf LFI}s presented in~\cite{marcos:00} (see also~\cite{Tax}). Moreover, it could be interesting to find interesting model-theoretic results for these logics, including \quci. For instance, the Keisler-Shelah Theorem, which was obtained by T. Ferguson for \qmbc\ in~\cite{fer:18}, should be valid in these stronger logics. It would be interesting to adapt the proof given  in~\cite{fer:18} to the semantical framework of triples.

\

\section*{Acknowledgments} 
Coniglio was financially supported by an individual research
grant from CNPq, Brazil (308524/2014-4).

\

\end{document}